\documentclass[12pt]{article}

\usepackage{amsmath,amsthm,amssymb}

\theoremstyle{plain}
\newtheorem{theorem}{Theorem}
\newtheorem{proposition}{Proposition}
\newtheorem{corollary}{Corollary}
\newtheorem{lemma}{Lemma}[section]
\theoremstyle{remark}
\newtheorem{rem}{\sc Remark}[section]

\makeatletter

\@addtoreset{equation}{section}
\makeatother

\newcommand{\al}{\alpha}
\newcommand{\de}{\delta}
\newcommand{\ep}{\varepsilon}

\newcommand{\om}{\omega}

\newcommand{\Ga}{\Gamma}
\newcommand{\La}{\Lambda}

\newcommand{\beq}{\begin{eqnarray*}}
\newcommand{\eeq}{\end{eqnarray*}}
\newcommand{\beqn}{\begin{equation}}
\newcommand{\eeqn}{\end{equation}}

\newcommand{\nin}{\noindent}
\newcommand{\pf}{\noindent {\it Proof. \,}}
\newcommand{\ds}{{\rm d}}
\newcommand{\as}{{\rm a}}

\newcommand{\vom}{\varOmega}
\def\v2{\vskip2mm}

\begin{document}

\begin{center}
{\Large  The two-sided exit problem \\  for a random walk on  $\mathbb{Z}$
 with infinite variance I }
\vskip4mm
{K\^ohei UCHIYAMA\footnote{Department of Mathematics,Tokyo Institute of Technology, Japan}} \\
\vskip2mm
\end{center}

\vskip6mm

\begin{abstract} Let $S=(S_n)$ be an oscillatory  random walk on the integer lattice  $\mathbb{Z}$ with i.i.d. increments. Let $V_{{\rm d}}(x)$ be the renewal function of the strictly descending ladder height process for $S$. We obtain several sufficient conditions---given in terms of the distribution function of the increment $S_1-S_0$---so that as $R\to\infty$
$$ (*) \quad P [  S\;  \mbox{leaves $[0,R]$ on its upper side}\, |\, S_0=x] \, \sim\,  V_{{\rm d}}(x)/V_{{\rm d}}(R)$$
uniformly for $0\leq x\leq R$. When  $S$ is attracted to a stable process of index $0<\al \leq 2$ and there exists $\rho= \lim P[S_n>0]$, the sufficient condition obtained are  also necessary  for $(*)$ and fulfilled  if and only if  $(\alpha\vee 1)\rho =1$, and some asymptotic estimates of  the probability on the left side of $(*)$  
are given  in case $(\alpha\vee 1)\rho \neq 1$. 
\end{abstract}



\vskip6mm
\section{Introduction}

Let $S= (S_n, n=0,1, 2, \ldots) $ be  a  random walk (r.w.)  on the integer lattice $\mathbb{Z}$ with i.i.d.  increments and  an initial point  $S_0$ which is an unspecified integer.   Let $X$ be a generic random variable having the same law as the increment $S_1-S_0$. 
For $x\in \mathbb{Z}$ let   $P_x$ denote  the law  of the r.w. $S$ started at $x$ and  $E_x$ the expectation by $P_x$;  the subscript $x$ is dropped from $P_x$ and $E_x$ if $x=0$. 
We suppose throughout the paper that $S$ is  irreducible (as a Markov chain on $\mathbb{Z}$) and {\it oscillating} (i.e., the sequence $(S_n)$ changes signs 
infinitely often with probability 1). 
For  a  subset $B\subset (-\infty,\infty)$ such that $B\cap \mathbb{Z} \neq \emptyset$,  denote by $\sigma_B$ (sometimes by  $\sigma B$ when it appears as a subscript) the first time when the r.w. $S$ visits  $B$ after time zero, namely $\sigma_B= \inf\{n\geq 1: S_n\in B\}$. We shall mainly work with the negative half line $(\infty,-1]$ which we denote by $\vom$:
$$\vom = (-\infty,-1].$$
Put
$$
T=\sigma_{\vom}.
$$
  This paper concerns the asymptotic form  as $R\to\infty$ of the probability of  
  $$\La_R =\{\sigma_{(R,\infty)}< T \},$$
the event that the r.w. $S$ exits from  an  interval  $[0,R]$ through the right boundary.
 ($R$ will always denote a positive integer).  The classical result given in \cite[Theorem 22.1]{S} says that if the variance $\sigma^2:= E X^2$ is finite, then $P_x(\La_R)-x/R \to 0$ uniformly for $0\leq x<R$. This can be refined to
\beqn\label{fun_res}
P_x(\La_R) \sim \frac{V_\ds(x)}{V_\ds(R)} \quad \mbox{uniformly for $0\leq x\leq R$ as} \; R\to\infty
\eeqn
as is given in \cite[Proposition 2.2]{U1dm} (the proof is easy: see  Section 3.1 of the present article). Here  $V_\ds(x)$ denotes the renewal function of the {\it weakly descending}  ladder height process for $S$ and $\sim$ means that the ratio of its two sides tends to unity. Our primary purpose  is to find a sufficient condition for (\ref{fun_res}) to hold that applies to a wide class of r.w.'s on $\mathbb{Z}$ with 
$\sigma^2= \infty$.

 If  the distribution of $ X$ is symmetric and belongs to the domain of normal attraction of a stable law with exponent $0<\alpha\leq 2$,  the problem is investigated by Kesten \cite{K}: for $1<\alpha\leq2$ he identifies the limit of $P_x[\sigma_{[R,\infty)}< T]$ as $x\wedge R\to\infty$ so that $x/R\to \lambda\in (0,1)$. Rogozin \cite{R2} studied the corresponding problem for strictly stable processes and derived an explicit analytic expression of the probability of the process exiting the unit interval on the right side (together with the law of the overshoot distribution) and,  as an application of the result, obtained the  scaled limit of $P_x(\La_R)$ for asymptotically stable r.w.'s (see Remark \ref{rem1.3}), extending the result of \cite{K} mentioned above. 
There are a few other results concerning the asymptotic behaviour of $P_x(\La_R)$. Griffin and MacConnel \cite{GM} gave a criterion for $P_x(\La_{2x})$  to approach unity (see Remark \ref{rem1}). 
 When $X$ is in the domain of attraction of a stable law, Doney \cite{Dsp} derived the asymptotic form of $P(\La_R)$.  Kesten and Maller \cite{KMdivg} gave an analytic condition for  $\liminf P_{x}(\La_{2x}) \wedge [1-P_{x}(\La_{2x})]>0$, i.e., for $S$ started at the centre of a long interval to exit it either side of it with positive probability.

Let $Z$ (resp. $\hat Z$) be  the first ladder height of the strictly ascending (resp. descending ladder) process: $Z= S_{\sigma[1,\infty)}$,  $\hat Z = S_T$, if $S_0=0$. Both $Z$ and $\hat Z$ are proper random variables since the r.w. $S$ is assumed to oscillate. The product  $E Z E|\hat Z|$ is finite  if $\sigma^2:=E  X^{2}<\infty$, and infinite otherwise (cf \cite[Section 17]{S}). 
$Z$ is said to be {\it relatively stable (r.s.)} if there is a norming constants  $B_n$ such that  $(Z_1+\cdots+Z_n)/B_n \to1$ in probability, where $Z_k$ are independent copies of $Z$. Let $X$ be a generic random variable having the same law as the increment $S_1-S_0$ and 
 $F$ be the distribution function of $X$.  If $E|X| <\infty$, put
$$\eta_+(x) = \int_x^\infty (1-F(t)) dt, \quad  m_+(x) = \int_0^x \eta_+(t)dt, \quad $$
for $x\geq 0$; define $\eta_-$ and $m_-$ similarly with $F(-t)$ in place of  $1-F(t)$
and let  $\eta = \eta_- +\eta_+$ and $m=m_-+m_+$. In the sequel we always assume 
$$\sigma^2=EX^2=\infty.$$
Because of the oscillation of $S$, $EX=0$ if $E|X|<\infty$. (Below, we shall  write $EX=0$ simply  with the understanding that $E|X|<\infty$ is implicitly assumed.) A positive Borel function $f$ on a neighbourhood of $\infty$ is said to be {\it slowly varying (s.v.)} (at infinity) if $\lim f(\lambda x)/f(x)=1$ for each  
$\lambda>1$. 
Our main result is stated in Theorem \ref{thm1}  below and some results complementary to it are provided in Proposition \ref{prop1.1}.

\vskip3mm\nin
\begin{theorem}\label{thm1} 
 Let $\sigma^2 =\infty$. In order for  (\ref{fun_res})  to be true each of the following are sufficient
\v2

{\rm (C1)}  \; $ E X=0$ and   ${\displaystyle m_+(x) /m(x)\to 0}$ as $x\to\infty$.
\v2

{\rm (C2)} \; $E X=0$ and $m$ is s.v. as $x\to\infty$.
\v2
{\rm (C3)} \;\; $Z$ is r.s. and. $V_\ds$ is s.v.
\v2
{\rm (C4)} \;\; $F(-x)/(1-F(x)\to0$ as $x\to\infty$  and $1-F(x)$ is  regularly varying  

 $\qquad\quad$ at infinity with index $-\al$, $0< \al <1$.
\end{theorem}
\v2

\v2
The sufficiency of (C1) is proved by the present author in \cite[Proposition 47]{Upot}.  The primary objective of this note is to verify the sufficiency of  (C2) to (C4);  proof will be made independently of the result under (C1).     

It is warned that the second condition in (C4)    does not mean that  the left tail   $F(-x)$ may  get small arbitrarily fast as $x\to\infty$,  
  the random walk $S$  being supposed to  oscillate  so that according to \cite{E}   the growth condition
$\int_{-\infty}^0 \big(|t| \big/\int_0^{|t|}(1-F(s)))ds\big)dF(t) = \infty$ must be satisfied (valid for   general r.w.s with $E[X; X>0] =\infty$). In case (C3) (\ref{fun_res}) is of interest only when $x/R\to0$, for if $V_\ds$ is s.v.,  (\ref{fun_res}) entails $P_x(\La_R)\to1$ whenever $x/R$ is bounded away from zero.

Condition (C2) holds if and only if $F$ belongs to the domain of attraction  of a centred normal law
\cite[Appendix A]{Upot}. If $F$ is attracted to a stable law of exponent $\al$, then (C3) hold if and only if $\al=1=\lim P[S_n>0]$ (see  \cite{Uexp}).

In what follows we shall omit  \lq$x\to \infty$' or \lq$R\to\infty$' when it is obvious. 




\begin{rem}\label{rem1}\,  Put 
$$A(x) = \int_0^x[1-F(t) -F(-t)]dt\quad  \mbox{and}\quad H(x)= P_x[|X|>x].$$
 Then condition (C3) holds if
\beqn\label {KMC}
\frac{A(x)}{xH(x)} \, \longrightarrow\, \infty.
\eeqn 
Indeed, (\ref{KMC}) is equivalent to the positive relative stability of $X$ 
according to  \cite[p.1478]{KM/ifp} and the latter implies that $Z$ is r.s. \cite[Theorem10]{R}. Moreover the condition
\beqn\label {KMC2}
A(x)/xF(-x) \, \longrightarrow\, \infty,
\eeqn 
apparently weaker than (\ref{KMC}), is equivalent to $\lim P_x(\La_{2x}) =1$  \cite[p.1477]{KM/ifp} and, therefore, entails that $V_\ds$ is s.v. (see (\ref{upp_bd} below). In view of Theorem \ref{thm1}, (C3) implies $P_x(\La_{2x}) \to 1$, hence (\ref{KMC2}). It remains unclear if (\ref{KMC}) holds under (C3).
\end{rem}


\begin{rem}\label{rem2}\, \, If $U_\ds$ denotes the renewal function of the strictly descending ladder height process, then $U_\ds(x)=V_\ds(x)/V_\ds(0)$. We write $v^0$ for $V_\ds(0)$, which is given by
$$1/v^0 =P[S_{\sigma(-\infty,0]}<0] = P[S_{\sigma[0,\infty)}>0]
= \exp\Big\{-\sum k^{-1}p^k(0)\Big\}$$
 (cf., e.g., \cite[Sections XII.9, XVIII.3]{F}).  Put  
$$\ell^*(x) = \int_0^x P[Z>t]dt \quad\mbox{and}\quad 
\hat \ell^*(x) = \frac1{v^0}\int_0^x P[-\hat Z>t]dt.$$
 Then
\beqn\label{asmp_U}
V_\ds(x) \sim
 \left\{\begin{array}{ll} 
a(x)\ell^*(x)  &\mbox{ in case \;  (C1)},\\[1mm]
x/\hat\ell^*(x)  &\mbox{ in case \;  (C2)}, \\[1mm]
{\displaystyle  \bigg(\al\int_x^\infty \frac{U_\as(t)F(-t)}{t}dt \bigg)^{-1}}  \quad &\mbox{ in cases \;  (C3), (C4),} 
\end{array}\right.
\eeqn 
and  in case (C4) (as for (C3)) $V_\ds$ is s.v.  so that $P_{R/2}(\La_R)\to1$.
Here  $\al =1$ in case (C3),  $U_\as  (x)$ denotes the renewal function of the strictly ascending ladder height process,  and 
$$a(x) =\sum_{n=0}^\infty \big(P[S_n=0] - P[S_n =- x]\big),$$
 the potential function of $S$  (if $S$ is recurrent).  [See \cite[Section 7.3]{Upot}, \cite{R} and Lemma \ref{lem3.1} for the first, second and third formulae of (\ref{asmp_U}), respectively.] 

 In each case of (\ref{asmp_U}), $\ell^*$ or $\hat \ell^*$ or $F(-\,\cdot)U_\as   $ is s.v.    
Under (C1)   $a(x)\asymp x/m(x)$  (cf. \cite{Upot}), where $\asymp$ means that  the ratio of its two sides is bounded away from zero and infinity. 
 If (\ref{KMC}) holds, the first formula of (\ref{asmp_U}) has the following extension (see Remark \ref{rem4.2}): if $p:= \lim [1-F(x)]/H(x), q=1-p$, $F$ is recurrent and $p\neq 1/2$ (in addition to (\ref{KMC})), then
$$V_\ds(x)\sim (1-q^{-1}p)a(x)\ell^*(x).$$
\end{rem}

\v2  
\v2
In all the cases (C1) to (C4) either  $Z$ is r.s.  or $P[S_n>0]\to 1$
(the latter is equivalent to  $P_x(\La_{2x})\to1$ under the present setting (cf. \cite{KM/ifp}; see \cite[Theorem 8(i)]{Upot}) for the case (C1)).  
Below we suppose the 
  asymptotic stability condition
 \[
 {\rm (AS)} 
 \left\{
 \begin{array}{ll}
 {\rm a.} \; X\; \mbox{is attracted to a stable law of exponent} \; 0<\alpha \leq 2.\\
   {\rm b.} \; EX=0\; \mbox{if}\; E|X|<\infty. \\
   {\rm c.} \mbox{ \;\;  there exists} \quad 
   \rho :=\lim P[S_n>0],
   \end{array}\right.
   \]
   and derive some upper and lower bounds of $P_x( \La_R)$.
Note that the assumption of   $S$ being oscillating is always in force.

Suppose condition (ASab)---the conjunction of (ASa) and (ASb)---to hold with $\al<2$. It then follows that 
 \beqn\label{stable}
 1-F(x) \sim pL(x)x^{-\alpha}  \quad \mbox{and} \quad F(-x) \sim q L(x)x^{-\alpha}
 \eeqn
 for some s.v. function $L$ and  constant $p=1-q\in [0,1]$.  Let the positive numbers $c_n$ be chosen so that $\lim H(c_n)n =(2-\alpha)/\alpha $. For $\al\neq1$, (ASc) is valid with $[1-\al^{-1}]\vee0 \leq \rho\leq \al^{-1}\wedge1$ and $S_n/c_n$ converges in law to a no-degenerate stable variable.
 For $\al=1$,  (ASc) holds with $\rho=1$ at least  if either $p<1/2, EX=0$ or $p>1/2, E|X|=\infty$; if $p=1/2$,  in order that  $S_n/c_n$ converges in law to a non-degenerate variable it is necessary and sufficient that  (ASc) holds with $0<\rho<1$ and this is equivalent to 
  the existence of  $\lim n E[\sin X/c_n] \in \mathbb{R}$. (Cf. \cite[Section XVII.5]{F}.) 
 
If  (AS) holds, then it follows that  $\al\rho\leq 1$ if $\al=2$;  $Z$ is  r.s. if and only if $\al\rho=1$;
for $\al=2$, $\rho=1/2$; for $1<\al<2$, (C1) $\,\Leftrightarrow\, p=0 \,\Leftrightarrow\,\al =\rho=1$, and that
  \beqn\label{EQS}
 {\rm (C2)} \,\Leftrightarrow\, \al=2;\;\; {\rm (C3)} \,\Leftrightarrow\,\al=\rho=1; \;\; {\rm (C4)} \,\Leftrightarrow\, \al<1=p \,\Rightarrow\, \rho=1;
  \eeqn
in particular, either one of (C1) to (C4) holds if and only if $(\al\vee1)\rho=1$.  We consider the case
$(\al\vee1)\rho<1$ in the next proposition which combined with Theorem \ref{thm1} shows that
 (\ref{fun_res}) holds   if and only if  $(\al\vee1)\rho=1$.   
\vskip3mm
\begin{proposition}\label{prop1.1} \,  Suppose (AS) to hold with $0< \al <2$.   
 Let $\de$  be a constant  arbitrarily chosen so that $ \frac12<\de <1$.
\v2
 {\rm (i)}\; If  $0<(\al\vee 1)\rho<1$ (equivalently,   $p>0$ for $\al>1$ and $ 0<\rho <1$ for $\al\leq 1$),   then there exist  constants  
 $\theta_*>0$ and $\theta^* <1$ such that for all sufficiently large $R$,
\beqn\label{th/th}
\theta_* \leq \frac{P_x(\La_R) V_\ds(R)}{V_\ds(x)}\leq \theta^*   \qquad \mbox{ for} \quad  0\leq x\leq \de R.
\eeqn

\v2

{\rm (ii)}\; Let $\al \leq 1$. In (\ref{th/th}), one can choose $\theta_*$ and $\theta^*$ so that $\theta^* \to 0$ as $\rho\to 0$ and $\theta_*\to 1$ as $\rho\to1$. \big[The statement \lq\lq $\,\theta^*\to 0$ as $\rho\to 0$'' means that for any $\ep>0$ there exists $\de >0$ such that $\theta^* <\ep$ for any \lq admissible' $F$ with $0<\rho<\de$, and similarly for  \lq\lq $\,\theta^*\to 1$ as $\rho\to1$''.\big] If $\rho=0$, then $P_x(\La_R)V_\ds(R)/V_\ds(x) \to 0$ uniformly for $0\leq x <\de R$, and $V_\ds(x)\sim x^\al/\hat\ell(x)$ for some s.v. function $\hat\ell$. [For $\al=1$, one can take $\hat\ell^*$ for $\hat\ell$ (see (\ref{5.30}).] 
\v2
\end{proposition}

\begin{rem}\label{rem1.3} Suppose that (AS) holds and that $0<\rho<1$ if $\alpha\neq1$ (i.e., limiting
stable law is strictly stable). Let $Y=(Y_t)_{t\geq 0}$ be the limiting stable process (which we suppose to be defined on the same probability space as $S$), denote by  $\La_r^Y$, $r>0$  the corresponding event for $Y$:
$$\La_r^Y = \big\{\mbox{the first exit  of $Y$ from the interval $[0,r]$ is on the right side} \big\},$$
and put  $Q(\xi) = P(\La_1^Y |\, Y_0=\xi)$ 
$(0< \xi<1)$.   Rogozin \cite{R2} established the overshoot law given in (\ref{Rog1}) below and from it  derived an analytic expression of $Q(\xi)$ which, if $0<\rho <1$,  must read 
\beqn\label{1.10}
Q(\xi) =\frac{1}{B(\alpha\rho,\alpha\hat\rho)}  \int_0^{\xi}t^{\al\hat\rho-1}(1-t)^{\al\rho-1}dt \quad (0<\xi<1),
\eeqn
where $\hat\rho =1-\rho$ and $B(s,t)=\Gamma(s)\Gamma(t)/\Gamma(s+t)$,  the beta function.  (For $\alpha<1$, $Q(\xi)=1$ or $0$ ($\xi>0$) according as $\rho=1$ or $0$.) By the functional limit theorem one deduces that $P_{\xi c_n}(\La_{c_n}) \to P[\La_1^Y\,\,|\, Y_0=\xi]$ for each $0<\xi<1$.  This convergence, though uniform since $P_x(\La_R)$ is monotone in  $x$ and $Q(\xi)$ is continuous, does not imply (\ref{th/th}). 

 Let $Z^Y\!(r)$ denote the overshoot of $Y$ as $Y$ crosses the level  $r$ from the left to the right of it (the analogue of $Z(R)$ given at the beginning of the next section).
 Then the overshoot law of Rogozin \cite{R2} mentioned above may read as follows:  if $0<\alpha\rho \leq 1$, then 
\beqn\label{Rog1}
P\big[Z^Y(1) >\eta,  \La^Y_1 \big|\,  Y_0=\xi\big] = \frac{\sin  \alpha\rho\pi}{\pi} (1-\xi)^{\alpha\rho}\xi^{\alpha\hat\rho} \int_\eta^\infty \frac{t^{-\alpha\rho}(t+1)^{-\alpha\hat \rho}}{t+1-\xi}dt
\eeqn
for $0<\xi <1$,  $\eta>0$. [Note that if $\alpha\rho=1$, the the RHS vanishes,  conforming to the fact that $Z$ is r.s, also that  the case $\alpha=1$ with $\rho\hat\rho=0$ is excluded.] 

 The verification of (\ref{1.10}), given in \cite{R2},  does not seem to be finished and the formula obtained there differs from (\ref{1.10})---$\rho$ and $\hat\rho$ are interchanged. Here we provide proof of (\ref{1.10}),  based on (\ref{Rog1}). Suggested by proof of Lemma 1 of Kesten \cite{K} we prove instead of 
  (\ref{1.10}) 
 \beqn\label{Kder}
 P(\La_{1+c}^Y\,|\, Y_0=1) = \frac{1}{B(\alpha\rho,\alpha\hat\rho)} \int_0^{1/(1+c)}t^{\al\hat\rho-1}(1-t)^{\al\rho-1}dt \quad\;\;  (c>0),
\eeqn
which is equivalent to (\ref{1.10}).
By  (\ref{Rog1}) the conditional probability on the LHS is expressed as 
$$ \frac{\sin \pi \alpha\rho}{\pi} c^{\alpha\rho} \int_0^\infty \frac{t^{-\alpha\rho}(t+1+c)^{-\alpha\hat \rho}}{t+c}dt.
$$
On changing the variable  $t=c/s$  this yields
$$ P(\La_{1+c}^Y\,|\, Y_0=1)= \frac{\sin \pi \alpha\rho}{\pi} \int_0^\infty \frac{s^{\alpha -1}(s+(1+s)c)^{-\alpha\hat \rho}}{1+s}ds.
$$
Using the identity $\int_0^\infty s^{\gamma-1}(1+s)^{-\nu}ds =B(\gamma,\nu-\gamma)$
 ($\nu >\gamma >0$) one can easily deduce that the above expression and the RHS of  (\ref{Kder}) equal 1 for $c=0$ and  their derivatives with respect to  $c$ coincide. 
  \end{rem}
 \v2
\begin{rem}\label{rem1.4}  With the help of  our Theorem \ref{thm1} we obtain from (\ref{Rog1})   that
there exists
$$Q^\circ(\xi,\eta) :=\lim_{R\to\infty} P_{\lfloor\xi R\rfloor}[Z(R) > \eta R, \La_R] \quad (0\leq \xi\leq 1, \eta >0),$$
for which 
$$Q^\circ(\xi,\eta) = P\big[Z^Y(1) >\eta,  \La^Y_1 \big|\,  Y_0=\xi\big] \quad \; \mbox{unless}\quad 1-\alpha=0=  \rho\hat\rho;\mbox{and}$$
$$Q^\circ(\xi,\eta) =
\left\{
\begin{array}{ll}
0       \qquad   \qquad   \qquad \;  \quad(\xi<1, \eta >0)    \quad &\mbox{if}\quad  \alpha\rho\in \{0,1\},\\
C(\xi)\eta^{-\alpha}\{1+o(1)\} \quad\mbox{as}\;\; \eta\to\infty    \quad &\mbox{if}\quad 0< \alpha \rho<1,
\end{array}\right.
$$
where $C(\xi)$ is a continuous function positive for  $0< \xi<1$ and $C(0)=C(1)= 0$. 
\end{rem} 
\v2\v2 
 
As mentioned above, the formula (\ref{1.10}) yields the explicit asymptotic form of $P_x(\La_R)$ if $x/R$ is bounded away from zero unless $\alpha=1$ and $\rho\hat\rho=0$, but  when  $x/R\to0$ it says only $P_x(\La_R) \to 0$.  One may especially be interested in the case $x=0$.   
  Doney \cite[Corollary 3]{Dsp} 
obtained under (AS) with $0<\rho<1$  a functional limit theorem of the scaled r.w. conditioned on various events (that depend on the behaviour of the r.w. only up  to $T$), and derived, as an application of it, that 
\beqn\label{BD2} 
P\big[ \max_{n<T} S_n \geq b_n\big] = P(\La_{b_n}) \sim c/n
\eeqn
where $c$ is a positive constant and $b_1, b_2,\ldots$ are  norming constants determined by the relation $nP[-\hat Z>b_n]\to 1$ if $\alpha\hat\rho<1$ and $\int_0^{b_n}P[-\hat Z>t]dt/b_n  \sim b_n/n$ if $\alpha\hat\rho=1$.  If (\ref{fun_res}) holds, taking  $x=0$ therein  gives  $P(\La_R)\sim v^0/V_\ds(R)$.
[It holds that $V_\ds(b_n)/n$ tends to a positive constant (see (\ref{U/Z}) of the next section) and one can easily ascertain that the equivalence in (\ref{BD2}) may be written as $P(\La_R)\sim c/V_\ds(R)$.]  

The corollary below gives the corresponding consequence of Proposition \ref{prop1.1} 
which improves on the result of \cite{Dsp} mentioned right above by extending it to the case $\rho \in \{0,1\}$ and providing information on how the limit depends on $\rho$. 
 
Suppose (AS) to hold.    We distinguish  the following three cases:
\v2
Case $I$: \quad  \,    $(\al\vee1) \rho =1$.
\v2
 Case $I\!I$: \quad \, $0< (\al\vee1) \rho <1$. 
 \v2
 Case $I\!I\!I$: \quad  $\rho=0$ 
 \v2\nin
 If $\al=2$, $I$ is always the case.  For $1<\al <2$, $I$ or $I\!I$ is the case according as $p=0$ or $p>0$. If $\al=1$, we have Case $I$ only if $p=1/2$, otherwise all the cases $I$ to  $I\!I\!I$ are possible. For $\al<1$, $I$ or $I\!I\!I$ is the case according as $p=1$ or $p=0$, otherwise we have Case $I\!I$.
 
The following result is obtained as a corollary of Theorem \ref{thm1} and Proposition \ref{prop1.1} except for the case  $\alpha>1$ of its last assertion.

\v2
 \begin{corollary}\label{cor1.1} Suppose (AS) to hold.   Then  
\vskip1mm
{\rm  (i)} there exists  $c =\lim P(\La_R)V_\ds(R)/v^0;$

\vskip1mm
{\rm  (ii)}    $c=1$  in Case $I$,  $0<c<1$ in Case $I\!I$ and $c=0$ in Case $I\!I\!I$;  

\vskip1mm
{\rm  (iii)}  $c\to0$ as  $\rho\to 0$ (necessarily $\alpha \leq 1$),  and \; $c\to 1$\,  as\,  $\rho\to1/(\al\vee 1)$. 
\end{corollary}

\v2\nin
\pf  As to the proof of (i) and (ii) Cases  $I$ and $I\!I\!I$ are covered by Theorem \ref{thm1} and by the second half of Proposition \ref{prop1.1}(ii), respectively,  and Case $I\!I$   by
the above mentioned  result of \cite{Dsp}  combined with Proposition \ref{prop1.1}.
The assertion (iii) follows from Proposition \ref{prop1.1}(ii) if $\al \leq 1$. In case  $1<\al <2$, we need to use a result from \cite{Uexit2}.  There we obtain the asymptotic estimates of 
$\lim_{R\to\infty} P_x[\sigma_y<T\,|\,\La_R]$ 
for $x<y$ in Proposition 1.4(i) of \cite{Uexit2} (deduced from  the estimates of $u_\as$ and $v_\ds$ gotten independently of the preset work), which entails that as $R\to\infty$ 
$$\qquad \frac{P(\La_R)V_\ds(R)}{v^0} \geq \frac{P[\sigma_{\{ R\}}<T]V_\ds(R)}{v^0}\,\longrightarrow\, \frac{\al-1}{\al(1-\rho)}.\qquad \qed
$$

 \begin{rem}\label{rem1.5}  (a) In \cite{Dsp} $EX=0$ is assumed if $\al=1$. However, by examining the arguments therein, one sees that this assumption is superfluous (not only for  Corollary 3 of \cite{Dsp} but for the main theorem of it), so that (\ref{BD2}) is valid under (AS) with $0<\rho<1$.
 
  (b)  In Corollary 3 of \cite{Dsp} some related results are stated, which are paraphrased as follows:
 $P(\La_R) \sim c/R$ with $c>0$ if and only if  $E|\hat Z|<\infty$; and $\int_0^x P(\La_t)dt$ is s.v. if
 $\al\hat\rho=1$ and $0<\al<2$. These follow immediately from our Corollary \ref{cor1.1} (recall that $V_\ds(x)\sim x/\hat\ell^*(x)$, $\lim \hat\ell^*(x) = E|\hat Z|\leq \infty$ and $\ell^*$ is s.v. if $\al\rho=1$).  
 \end{rem}

 \v2
 The present paper has been written simultaneously with \cite{Upot}. The two are closely related, and we make use of some results from one in the other: we shall apply Proposition 40 of \cite{Upot} in the proof of Theorem \ref{thm1} (in case (C2)) of the present one, while Lemmas \ref{lem4.4} and \ref{lem6.1}  of the present one are used in Section 7.4 of \cite{Upot}. These however do not give rise to any risk of  circular arguments.


The rest of this paper is organised as follows. In Section 2, we state several  known facts and provide a lemma;  these will be fundamental in the later discussions. Proof of  Theorem \ref{thm1} is given in Sections 3.  In Section 4 we assume (C3) or (C4) to hold and prove lemmas as to miscellaneous matters.  In Section 5 we prove Proposition \ref{prop1.1} by further showing several lemmas.

\section{Preliminaries}

Let $v_\ds$ denote the renewal sequence for the weakly descending ladder height process:  $v_\ds(x) = V_\ds(x)- V_\ds(x-1)$ and $v_\ds(0)=v^0$. Here and in the sequel $V_\ds(x)$ is set to be zero for $x\leq -1$. Similarly let $U_\as   (x)$ and $u_\as   (x)$ be the renewal function and sequence for the {\it strictly ascending} ladder height process. 
We shall be also concerned with  the overshoot which we define by  
$$Z(R) = S_{\sigma[R+1,\infty)} -R.$$
 
The function $V_\ds$ is harmonic for the r.w. $S$  killed as it enters $\vom$ in the sense that
\beqn\label{Hrm/V} 
\sum_{y=0}^\infty V_\ds(y)p(y-x) = V_\ds(x) \qquad (x\geq 0),
\eeqn
where $p(x) = P[X=x]$ \cite[Proposition 19.5]{S}, so that the process $M_n= V_\ds(S_{n\wedge T})$ 
 is a martingale under  $P_x$ for  $x\geq 0$.  By  the optional stopping theorem  one 
deduces  
\beqn\label{mart}
E_x[ V_\ds(S_{\sigma(R,\infty)});  \La_R]  = V_\ds(x).
\eeqn
Indeed, on passing to the limit in  $E_x[M_{n\wedge \sigma(R,\infty)}]= V_\ds(x)$  Fatou's lemma  shows that  the expectation on the LHS is less than or equal to  $V_\ds(x)$, which fact  in turn shows that  the martingale $M_n$ is uniformly integrable since $M_{n\wedge \sigma_{(R,\infty)}}$ is bounded
by the summable random variable $V_\ds(R)\vee M_{ \sigma_{(R,\infty)}}$. 
Obviously  the expectation in (\ref{mart}) is not less than   $V_\ds(R)P_x(\La_R) $ so that
\beqn\label{upp_bd}
P_x(\La_R) \leq V_\ds(x)/V_\ds(R).
\eeqn 
If either  $V_\ds$ or $\int_0^x P[-\hat Z>t]dt$ is regularly varying,  it follows \cite[Eq(8.6.6)]{BGT}  that
 \beqn\label{U/Z}
 \frac{V_\ds(x)}{x v^0}\int_0^xP[-\hat Z>t] dt\, \longrightarrow\, \frac1{\Ga(1+\alpha\hat\rho)\Ga(2-\alpha\hat\rho)}.
 \eeqn

For a non-empty $B\subset \mathbb{Z}$ we define  the Green function  $g_B(x,y)$ of the r.w. killed as it hits $B$ by
\beqn\label{GF}
g_B(x,y) = \sum_{n=0}^\infty P_x[ S_n=y, n<\sigma_B].
\eeqn
(Thus if $x\in B$,  $g_B(x,y)$ is equal to $\de_{x,y}$ for $y\in B$ and to $ \sum p(z-x) g_B(z,y)$ for $y\notin B$.)  We shall repeatedly apply Spitzer's formula
\beqn\label{G-ft}
g_{\vom}(x,y) =\sum_{k=0}^{x\wedge y} v_\ds(x-k)u_\as   (y-k)\quad \mbox{for} \; x,y \geq 0
\eeqn
\cite[Propositions 18.7,  19.3]{S}.  It follows that for $x\geq 1$, 
\beqn\label{hat-Z}
P[\hat Z = -x] = \sum_{y=0}^\infty g_{\vom}(0,y)p(-x-y) = v^0  \sum_{y=0}^\infty u_\as   (y)p(-x-y)
\eeqn
and, by the duality relation $g_{[1,\infty)}(-x,-y) = g_{\vom}(y,x)$,
\beqn\label{Z}
P[Z = x] = \sum_{y=0}^\infty g_{[1,\infty)}(0,-y)p(x+y) = \sum_{y=0}^\infty v_\ds(y)p(x+y).
\eeqn
(See \cite[Eq(XII.3.6a)]{F} for a direct derivation.)

Let $\ell^*$ be as in Remark \ref{rem2}, namely
$ \ell^*(x) =\int_0^x P[Z>t]dt.$ It is known that  $Z$ is r.s. if and only if $\ell^*$ is s.v. \cite{R} and if this is the case 
\beqn\label{ell*} u_\as   (x) \sim 1/\ell^*(x) \eeqn
\cite[Appendix(B)]{Upot}, \cite{Urenwl}.  

\begin{lemma}\label{lem1}  For $0\leq x\leq R$, 
$$\sum_{y=0}^R g_{\vom}(x,y) \leq  V_\ds(x)U_\as   (R).$$
Take a positive constant  $\de<1$  arbitrarily.   If $U_\as   $ is regularly varying with positive index, then  there exists a constant $c>0$ such that  for $0\leq x\leq \de R$
$$\sum_{y=0}^R g_{\vom}(x,y) \geq c\,  V_\ds(x)U_\as   (R),$$
and if $U_\as$ is s.v., then uniformly for $0\leq  x<\de R$
$$\sum_{y=0}^R g_{\vom}(x,y) \sim V_\ds(x)U_\as   (R).$$
\end{lemma}

 \pf \, Split the sum on the LHS of the inequality of the lemma at $x$. By  (\ref{G-ft}) the sum over $y > x$  and that over $y\leq x$ are equal to, respectively,
  $$\sum_{y=x+1}^R\sum_{k=0}^x v_\ds(k) u_\as   (y-x+k) = \sum_{k=0}^x v_\ds(k)\big[U_\as   (R-x+k)- U_\as   (k) \big]  $$
  and
   $$\sum_{k=0}^x\sum_{y=k}^x v_\ds(x-y+k) u_\as   (k) = \sum_{k=0}^x \big[V_\ds(x) -V_\ds(k-1)\big]u_\as   (k).$$
 Summing by parts yields
    $ \sum_{k=0}^x \big[V_\ds(x) -V_\ds(k-1)\big]u_\as   (k) =\sum_{k=0}^x v_\ds(k) U_\as   (k) ,$
    so that
$$\sum_{y=0}^R g_{\vom}(x,y) =\sum_{k=0}^x v_\ds(k) U_\as   (R-x+k) $$
by which the assertions of the lemma follow immediately.
 \qed

\section{Proof of Theorem \ref{thm1}}

Proof of Theorem \ref{thm1} is given  in two subsections, the first  for the case (C2) and the second for cases (C3) and (C4).   
\v2
\v2
{\bf 3.1. Case (C2).} 
\v2
Suppose $EX=0$ and recall $\eta_\pm(x) =\int_x^\infty P[\pm X>t]dt$ and $\eta= \eta_-+\eta_+$.
It follows that if  $\lim x\eta_+(x)/m(x) = 0$, then $Z$ is r.s. \cite[Proposition 40]{Upot}, so that $u_\as   (x)\sim 1/\ell^*(x)$ is s.v.   owing to (\ref{ell*}). 
Thus (C2), which is equivalent to $\lim x\eta(x)/m(x) = 0$,  entails  both  $v_\ds$ and $u_\as$ are s.v.  so that   $\sum_{k=0}^y v_\ds(k)u_\as   (k) \sim yv_\ds(y)u_\as   (y) \sim V_\ds(y)u_\as   (y) $ ($y\to\infty$) and one can easily deduce that 
as $R\to\infty$
$${g_{\vom}(x,R)} \sim \sum_{k=0}^x \frac{v_\ds(k)}{\ell^*(R-x+k)}\sim   \frac{V_\ds(x)}{\ell^*(R)} \quad\mbox{uniformly for} \;\; 0\leq x<R.$$ 
This leads to the lower bound 
\beqn\label{-2}
P_x(\La_R) \geq P_x[\sigma_{\{R\}}<T] = \frac{g_{\vom}(x,R)} {g_{\vom}(R,R)}= \frac{V_\ds(x)}{V_\ds(R)}\{1+o(1)\},
\eeqn
which combined with (\ref{upp_bd}) shows (\ref{fun_res}) as desired.
\v2


\v2
\v2
{\bf 3.2.  Cases (C3) and (C4).}\; 
\v2
Each of conditions (C3) and (C4) implies that $U_\as   (x)/x^\al$ (necessarily $0  < \al \leq 1$) is s.v. and that $\lim P[S_n>0]=1$.
 Let $\ell$ be a s.v. function such that
 \beqn\label{-1}
 U_\as   (x) \sim
 x^\al/\ell(x) \quad \mbox{if} \;\; \al<1 \quad \mbox{and} \quad \mbox{if} \;\; \al=1, \quad \ell =\ell^*.
 \eeqn
Here and in the sequel   $\al=1$ in cases  (C3). We then bring in  the function
\beqn\label{3.0}
\ell_\sharp(t) = \al\int_t^\infty \frac{s^{\al-1}F(-s)}{\ell (s)}ds \quad (t >0),
\eeqn
where the summability of the integral is assured by $U_\as   (x)= \sum_{-\infty}^0 U_\as   (y)p(y+x)$, the dual relation of (\ref{Hrm/V}). 
In cases (C1)  and (C2) we have $P_x[\sigma_{\{R\}} <T] \sim V_\ds(x)/V_\ds(R)$ which combined with  (\ref{upp_bd}) entails  (\ref{fun_res}), but in cases (C3) and (C4) this
equivalence does not generally hold   (see Proposition 1.4 of \cite{Uexit2}) and for proof of (\ref{fun_res})  we shall  make use of (\ref{mart}) not only for the upper bound but for the lower bound.


\begin{lemma}\label{lem3.1}   Suppose that either (C3) or (C4) holds.  Then $\ell_\sharp(x)$ is s.v., and
$$(a)\;\; P[-\hat Z\geq x]  \sim v^0\ell_\sharp(x)\quad\mbox{and}\quad  (b)\;\; V_\ds(x)\sim 1/\ell_\sharp(x).
$$
\end{lemma}

\v2\nin
\pf  We have only to show (a), in addition to $\ell_\sharp$ being s.v.,  (a)  entailing (b) in view of (\ref{U/Z}) applied with $\alpha\hat\rho=0$. If (C3) or (C4) holds, then $V_\ds$ is s.v.\,(cf. \cite{Urenwl} in case (C4)), hence $P[-\hat Z\geq x]$ is also s.v. \cite{R}. Because of this fact together with the expression $P[-\hat Z\geq x] = v^0\sum_{y=0}^\infty u_\as    (y)F(-x-y)$, the assertion follows from a general result involving regularly varying functions that we state and verify in Lemma \ref{lem8.1} of Appendix (A) in order not to break the continuation of the discourse. [Under (C3), Lemma \ref{lem3.1} is a special case of Lemma 50 of \cite{Upot} where a different proof is found.]  \qed

 
  
 \v2
 
\begin{lemma}\label{lem3.2} {\rm (i)}   Let   (C3)   hold.    Then for any $\de>0$,
\beqn\label{3.1} P[Z>x] \geq  V_\ds(x)[1-F(x+\de x)]\{1+o(1)\}.
\eeqn
If \;  $\limsup [1-F(\lambda x)]/[1-F(x)] <1$ for some $\lambda>1$ in addition, then
\beqn\label{3.2}
P[Z>x] \leq  V_\ds(x)[1-F(x)]\{1+o(1)\}.
\eeqn

{\rm (ii)}\;  If (C4) holds, then \;\;  $P[Z>x] \sim V_\ds(x)[1-F(x)]$.
\end{lemma}
\pf  First we prove (ii).  Suppose (C4) to hold and let $1-F(x)\sim L_+(x)/x^\al$ with a s.v. function $L_+$. By (\ref{Z}) 
\beqn\label{3.3}
P[Z>x] = \sum_{y=0}^\infty v_\ds(y)\big[1-F(x+y)\big].
\eeqn
  For $M>1$,  the sum over $y\leq Mx$ is
asymptotically equivalent to   
$$ \sum_{y=0}^{Mx} \frac{v_\ds(y)L_+(x+y)}{(x+y)^\al}\sim 
L_+(x) \sum_{y=0}^{Mx} \frac{v_\ds(y)}{(x+y)^\al} =  \frac{V_\ds(x)L_+(x)}{x^\al}\big\{1 +  o(1)\big\}, $$
since $V_\ds$ is s.v.\,owing to Lemma \ref{lem3.1} and hence $\sum_{y=\ep x}^{Mx} v_\ds(y) = o(V_\ds(x))$ for any $\ep >0$.
Choosing  $L_+$  so that $L_+'(x)=o(L_+(x)/x)$ we  also see that the sum over $y>Mx$ is at most  $V_\ds(x)[1- F(Mx)]\{1+o(1)\}$.  Since $M$ can be arbitrarily large, this concludes the asserted equivalence. 

The inequality (\ref{3.1}) is obtained by restricting the range of summation in (\ref{3.1}) to $y\leq \de x$.  (\ref{3.2}) follows from
\beqn\label{3.4}
P[Z>x] \leq  V_\ds(x)[1-F(x)] + \sum_{y=x}^\infty v_\ds(y)\big[1-F(y)\big].
\eeqn
Indeed, under the assumption for (\ref{3.2}), the sum of the infinite series is $o\big(V_\ds(x)[1-F(x)]\big)$ according to a general result, Lemma \ref{lem8.2} of Appendix (A), concerning s.v. functions.
 \qed
\v2
\begin{rem}\label{rem3.1} Lemma \ref{lem3.2} says that if  the positive tail of $F$ satisfies a mild  regularity condition (such as (\ref{stable})), then  $P[Z>x] \sim [1-F(x)]/ \ell_\sharp(x)$ under (C3). As to the integrals  $\int_0^x [1-F(t)]/\ell_\sharp(t)dt$ and  $\int_0^x P[Z> t]\ell_\sharp(t)dt$, we have the corresponding relations with no additional condition assumed (see Lemmas \ref{lem3.5} and \ref{lem3.6}). 
\end{rem}
The slow variation of $\ell^*(x)=\int_0^xP[Z>t]dt$ implies  $P[Z>x]=o(\ell^*(x)/x)$ so that
\beqn\label{3.5}
 Z\; \mbox{is r.s.}\, \Longrightarrow\, \lim U_\as   (x)P[Z>x]=0.
\eeqn
\begin{lemma}\label{lem3.3} 
 {\rm (i)} \;Under (C3) 
 \beqn\label{VUH}
 V_\ds(x)U_\as   (x) H(x) \to 0.
 \eeqn
 \v2
 {\rm (ii)}\; Under (C4), \quad $V_\ds(x)U_\as   (x)[1-F(x)] \to (\sin \pi\al)/\pi\al$.
\end{lemma}
\v2
One may compare the results above with the known one under  (AS) with $0<\rho<1$ stated in
(\ref{nUV}) below;  (i) will be refined in Lemma \ref {lem6.1}(iii) under (AS).
\v2\nin
\pf\,  Under (C4),  Lemma \ref{lem3.2}(ii) says $P[Z>x] \sim V_\ds(x)[1-F(x)]$ so that by Lemma \ref{lem3.1} $P[Z>x]$ is regularly varying with index $-\al$, which implies that $P[Z>x]U_\as(x)$ approaches $(\sin \pi\al)/\pi\al$ (cf. \cite[Eq(8.6.4)]{BGT}. Hence we have (ii). 

 As for (i),  on recalling the definition of  $\ell_\sharp$, the slow variation of $\ell_\sharp$ entails  $xF(-x)/\ell^*(x) =o(\ell_\sharp(x))$, so that on the one hand
 \beqn\label{3.6}
 V_\ds(x)U_\as   (x)F(-x) \sim V_\ds(x)xF(-x)/\ell^*(x) \sim xF(-x)/[\ell^*(x)\ell_\sharp(x)] \to 0.
 \eeqn
On the other hand,  applying  (\ref{3.4})---with $\de=1$ and with $\frac12 x$ in place of $x$---we deduce
 \beq
   V_\ds(x)U_\as   (x)[1- F(x)] &\sim& U_\as   (x)V_\ds({\textstyle \frac12}x)[1-F(x)]\\ &\leq& U_\as   (x)P[Z>{\textstyle \frac12}x]\{1+o(1)\}  \to 0.
\eeq
where the convergence to zero is due to (\ref{3.5}). Thus  we have (i). \qed
\v2
\begin{rem}\label{rem3.2}
(\ref{VUH}) holds also under (C2). Indeed,  we know $V_\ds(x)U_\as(x) \sim x^2\big/\int_{0}^x tH(t)dt$ (see (\ref{nUV})), and hence (\ref{VUH}) follows  since $x^2 H(x) = o\big(\int^x_0tH(t)dt\big)$. 
Under (C1),  (\ref{VUH}) may fail to hold, but it holds  (see \cite[Proposition 50(i)]{Upot}) that 
$$V_\ds(x)U_\as   (x)[1-F(x)] \to 0.$$
\end{rem}


\begin{lemma}\label{lem3.4}   {\rm (i)} \;Under (C3), for each $\ep >0$ as $R\to\infty$ 
\beqn\label{3.7}
 E_x\big[ V_\ds(S_{\sigma(R,\infty)}); Z(R)\geq  \ep R,\, \La_R\big]/V_\ds(x) \,\longrightarrow\, 0\quad\mbox{uniformly for \;\; $0\leq x\leq R.$} \eeqn
 
 {\rm (ii)}\; \;Under (C4),  as $M\wedge R\to\infty$ 
\beqn\label{3.8}
 E_x\big[ V_\ds(S_{\sigma(R,\infty)}); Z(R)\geq  MR,\, \La_R\big]/V_\ds(x) \,\longrightarrow\, 0\quad\mbox{uniformly for \;\; $0\leq x\leq R.$} \eeqn
 \end{lemma}
 \v2
 \pf \, We first prove (ii).  The expectation on the LHS of (\ref{3.8}) is less than
   \beqn\label{3.9}
  \sum_{w\geq R+ MR} \sum_{z=0}^{R-1} g_{\vom}(x,z) p(w-z)V_\ds(w)
  \eeqn
because of the trivial inequality $g_{\mathbb{Z}\setminus [0,R)}(x,z) < g_{\vom}(x,z)$.
 By Lemma \ref{lem3.1}  $V_\ds(x) \sim 1/\ell_\sharp(x)$  and the derivative  $(1/\ell_\sharp)'(x)=o(1/x\ell_\sharp(x))$. Hence  
 for each  $M>0$ there exists a constant  $R_0$ such that for all $R>R_0$ and $z<R$,
 \begin{eqnarray}\label{3.10}
  \sum_{w\geq R+MR} p(w-z)V_\ds(w) 
&\leq & \frac{1- F(MR)}{ \ell_\sharp(R)\{1+o(1)\}} 
+  \sum_{w\geq MR} \frac{1-F(w)}{w\ell_\sharp(w)}\times o(1) \\
& = &\frac{1-F(R)}{M^\al  \ell_\sharp(R)}\{1+o(1)\}
 \leq 2 \frac{V_\ds(R) [1-F(R)]}{M^\al},  \nonumber
\end{eqnarray}
where we employ the summation by parts  for the first inequality and (C4) for the equality.
  Owing to Lemma \ref{lem3.3} the last member is at most $3/M^\al U_\as   (R)$  and  hence  on  applying  Lemma \ref{lem1}  the double sum in  (\ref{3.9}) is dominated by  $3 V_\ds(x)/M^\al$, showing (\ref{3.8}). 
  
  In case (C3), by the inequality (\ref{3.1})  the sum of the infinite series in the second member
  of (\ref{3.10}) is at most a constant multiple of 
  $\sum_{w\geq MR} P[Z\geq \frac12 w]/w$, which is $o(\ell^*(R)/R)$ as one can easily check by summing  by parts. This together with the equivalence $U_\as   (x)\sim x/\ell^*(x)$ and Lemma \ref{lem3.3}(i) shows that the sum on the left-most member of (\ref{3.10}) is $o(1/U_\as   (R))$ and on returning to (\ref{3.9}) we conclude (\ref{3.7}) by Lemma \ref{lem1} again. 
  \qed
 
 
\v2
{\bf Proof of Theorem \ref{thm1} in cases (C3), (C4).}     We have the upper bound (\ref{upp_bd}) of $P_x(\La_R)$. To obtain the lower bound we observe that for $x\leq R/2$,
 $$\frac{V_\ds(R)}{1+o(1)} \geq   E_x\big[ V_\ds(S_{\sigma(R,\infty)}); Z(R)< MR \,\big|\,  \La_R \big]
 \geq \frac{ (1-\ep)V_\ds(x)}{P_x(\La_R)\{1+o(1)\}}:$$
 the left-had inequality follows immediately from the slow variation of $V_\ds$ (see Lemma \ref{lem3.1} in case (C4)), and the right-hand one from (\ref{mart}) and Lemma \ref{lem3.4}. The inequalities above 
 yield the lower bound for $x \leq R/2$ since $\ep$ can be made arbitrarily small.
In view of the  slow variation of  $V_\ds$
the result for $\lfloor R/2\rfloor \leq x <R$ follows from  $P_{\lfloor R/2\rfloor}(\La_R)\to1$. 
 The proof is complete. \qed

\v2
 
 \begin{rem}\label{rem3.3}
 Theorem \ref{thm1} together with Lemma \ref{lem3.4}  shows that uniformly for $0\leq x\leq R$
\beqn\label{3.11}
\begin{array}{ll}
\mbox{if  (C3) holds,}\quad  P_x[Z(R) > \ep R\,|\, \La_R] \to 0 \quad &\mbox{as \; $R\to\infty$ \; for each \; $\ep >0$}; \\[4mm]
\mbox{if  (C4) holds,}  \quad  P_x[Z(R) > M R\,|\, \La_R] \to 0 \quad &\mbox{as \;\; $R\wedge M\to\infty$}. 
 \end{array}
 \eeqn
\end{rem}
\v2

\section{ Miscellaneous lemmas under (C3), (C4).}

In this section we suppose (C3) or (C4) to hold and provide several lemmas for later 
citations.  In proof of the following two lemmas we put for $x\geq0$
\beqn\label{1-F}
H_+(x) = 1-F(x).
\eeqn

\begin{lemma}\label{lem3.5} \, Under (C3), $\ell^*(x)\sim \int_0^x V_\ds(t) [1-F(t)] dt$.
\end{lemma}
\pf  Summation by parts yields
$$P[Z>s]  = \sum_{y=0}^\infty v_\ds(y)[1-F(y+s)] = \sum_{y=0}^\infty V_\ds(y)p(y+1+ \lfloor s\rfloor).
$$
Since $\int_0^x p(y+1+ \lfloor s\rfloor)ds = F(y+1+x) -F(y) +O(p(\lfloor x\rfloor +y+1))$, putting
$$V(x)=1/\ell_\sharp(x)\quad\mbox{and}\quad \tilde\ell(x) =\int_0^\infty V(t)[F(t+x) -F(t)]dt,
$$
we have $\ell^*(x) =\int_0^x P[Z>s]ds =\tilde\ell(x) + o(1)$.  By splitting the integral at $x$ and employing a change of variables a simple manipulation yields
\beqn\label{3.12}
\tilde\ell(x) =\int_0^x V(t)H_+(t)dt + \int_0^\infty \big(V(y+x)-V(y)\big) H_+(t+x)dt.
\eeqn
In a similar way, taking the difference for the integral defining $\tilde \ell(x)$, we deduce
\begin{eqnarray}\label{3.13}
\tilde\ell(x)- \tilde\ell({\textstyle \frac12 x}) &=& \int_0^\infty V(t)\big[H_+(t+{\textstyle \frac12 x})- H_+(t+x)\big]dt\\
&=& \int_0^{x/2} V(t)H_+(t+{\textstyle \frac12 x})dt + \int_0^\infty \big(V(t+{\textstyle \frac12 x})-V(t)\big) H_+(t+x)dt. \nonumber
\end{eqnarray}
Since the integrals on the last member of (\ref{3.13}) are positive and $\tilde \ell$ is s.v., both must be of smaller order than $\tilde\ell$, in particular the second integral restricted to $t\geq x$ is $o(\tilde\ell(x))$, which  fact leads to
$$\int_0^\infty \big(V(t+x)-V(t)\big) \mu_+(t+x)dt =o(\tilde\ell(x)),
$$
for by the slow variation of $V$ as well as  of $1/\ell^*$ and the monotonicity of  $\mu_-(t)$, it follows  that if $t\geq \frac12 x$, then 
$V(t+x)-V(t) = \int_t^{t+x} V^2(s)F(-s)ds/\ell^*(s)\sim \frac{V^2(t)}{\ell^*(t)}\int_0^x F(-s-t)ds$, which shows
$$V(t+x)-V(t)  \leq 2\big[V(t+{\textstyle \frac12 x)} -V(t)\big]\{1+o(1)\}.
$$
Returning to (\ref{3.12}) we can now conclude $\tilde\ell(x)\sim\int_0^xV(t)\mu_+(t)dt$, hence by
Lemma \ref{lem3.1} $\ell^*(x)\sim\int_0^xV_\ds(t)\mu_+(t) dt$
as desired. \qed

\begin{lemma}\label{lem3.6} Suppose (C3) to hold. Then
\beqn\label{3.14}
\int_0^x \ell_\sharp(t) P[Z>t]dt \sim \int_0^x [1-F(t)]dt;
\eeqn
and if either $E[X; X>0]$ or $E[|X|; X<0]$ is finite in addition, then 
$$E|X|<\infty, \quad \int_0^\infty \ell_\sharp(t)P[Z>t]dt = E[X;X>0]$$
 and
\beqn\label{3.15}
\int_x^\infty \ell_\sharp(t)P[Z>t]dt \sim \int_x^\infty [1-F(t)]dt.
\eeqn
\end{lemma}
\pf  Denote by $J(x)$ the integral on the LHS of (\ref{3.14}). We first derive the lower bound for it. By  (\ref{3.1}) and $\lim \ell_\sharp(x)V_\ds(x)\to 1$, for any $\lambda >1$,
$$\frac{J(x)}{1+o(1)} \geq \int_0^x \ell_\sharp(t) V_\ds(t)H_+(\lambda t)dt =\frac1{\lambda}\int_0^{\lambda t}H_+(t)dt\{1+o(1)\}.
$$
Since $\lambda$ may be arbitrarily close to unity, we have $J(x)\geq \int_0^x H_+(t)dt\{1+o(1)\}$.

According to Lemma \ref{lem8.3} given in Appendix (A), by the slow variation of $V_\ds$ and the monotonicity of $H_+$ it follows that
$$\sum_{y=x}^\infty v_\ds(y)H_+(y) = o\bigg(V_\ds(x)H_+(x) + \int_x^\infty V_\ds(t)H_+(t) \frac{dt}{t}\bigg).
$$
In view of (\ref{3.4}), it therefore suffices for the upper bound to show that
\beqn\label{3.16}
\int_0^x \ell_\sharp(t)dt\int_t^\infty V_\ds(t)H_+(t)\frac{dt}{t} \leq C\int_0^x H_+(t)dt.
\eeqn
Integrating by parts we may write the LHS as
$$\bigg[x\ell_\sharp(x) \int_t^\infty V_\ds(t)H_+(t)\frac{dt}{t}+\int_0^x H_+(t)dt\bigg]\{1+o(1)\}.
$$
Integrating by parts again and using Lemma \ref{lem3.5} one can easily deduce that the first integral in
the big square brackets is less than $\int_x^\infty\ell^*(t)t^{-2}\{1+o(1)\}dt \sim \ell^*(x)/x\leq V_\ds(x)\int_0^x H_+(t)dt\times O(1/x)$. Thus we obtain (\ref{3.16}), hence (\ref{3.14}).

 Let $E[X; X>0] \wedge E[|X|; X<0] <\infty$. One can show (\ref{3.15}) in an analogous way to the above, but, under $E|X|<\infty$--- so that $EX=0$, (\ref{3.15}) follows immediately from (\ref{3.14}) since
 $$\int_0^\infty \ell_\sharp(t)P[Z>t]dt =\int_0^\infty \frac{F(-s)}{\ell^*(s)}ds\int_0^s P[Z>t]dt = \int_0^\infty F(-s)ds = \int_0^\infty H_+(t)dt.$$
 The verification of $E|X|<\infty$ also follows from (\ref{3.14}) by using  an estimate of  $A(x)$ and
 is postponed to the proof of the next lemma. \qed

  Let $\ell_\sharp$ be  the function defined by (\ref{3.0}) in the preceding section.
\begin{lemma}\label{lem4.4}  \, Suppose that either (C3) or (C4) holds. Then
 \beqn\label{id_A}
\ell^*(t) \ell_\sharp(t)= -\int_0^t F(-s)ds +\int_0^t P[Z>s]\ell_\sharp(s)ds
=A(t) + o\bigg(\int_0^t [1-F(s)]ds \bigg)
 \eeqn
and in case $EX=0$,  both $\eta_-$ and $\eta$ are s.v. and
\beqn\label{id_A2}
\ell^*(t) \ell_\sharp(t) = \int_t^\infty\big[F(-s) - P[Z>s]\ell_\sharp(s)\big]ds=A(t) +o(\eta_+(t)).\eeqn
 \end{lemma}
\pf
Recall  $\ell^*(t)=\int_0^tP[Z>s]ds$  and note  that as $t\downarrow 0$, $\ell^*(t) =O(t)$ and  $\ell_\sharp(t) \sim O(\log\, 1/t)$. Then  integration by parts leads to $\int_0^tP[Z>s]\ell_\sharp(s)ds =\ell^*(t)\ell_\sharp(t)+ \int_0^t F(-s)ds$, hence the first equality of (\ref{id_A}). The second one follows from (\ref{3.14}). [Note that   in case (C4) (\ref{3.14}) is immediate from Lemma \ref{lem3.2}(ii).

Here we complete  the proof of the second half of Lemma \ref{lem3.6}. Indeed, 
if $E[|X|;X<0]<\infty$, then (\ref{id_A}) verifies $E[X;X>0]<\infty$ as is clear from 
\beqn\label{ell-ell<}
 \ell^*(t)\ell_\sharp(t)< \int_t^\infty F(-s)ds\, (\leq \infty),
 \eeqn
 while  if   $E[X;X>0]<\infty$, then $E[|X|;X<0]<\infty$, since the right-most member of (\ref{id_A}) is positive.

Let $EX=0$. Then,  because of (\ref{ell-ell<}), the first equality of (\ref{id_A2}) is deduced from the identity
\beqn\label{el_el}
[\ell^*(t) \ell_\sharp(t)]'= P[Z>t]\ell_\sharp(t)  -  F(-t)
\eeqn
(valid almost every $t>1$). For the second one we have only to apply (\ref{3.15}).  Lemma \ref{lem3.3}(i) together with (b) of Lemma \ref{lem3.1} entails that under (C3)
$$xP[|X|\geq x] <\!< \ell^*(x)\ell_\sharp(x)\leq \eta_-(x).$$
 Thus both $\eta_-$ and $\eta$ are $s.v.$ \qed
 
 \v2
 \begin{rem}\label{rem4.2}  Suppose that (\ref{KMC}) and  $[1-F(x)]/H(x) \to p$,  and let $q=1-p$.   Then 
 combining Lemmas \ref{lem3.1} and \ref{lem4.4}  shows that for  $p\neq1/2$,  
 $$A(x)\sim \ell_\sharp(x)\ell^*(x) \sim x/[U_\as(x)V_\ds(x)].$$
 If  $F$ is recurrent (necessarily $p\leq 1/2$), we have 
 $a(x) \sim \int_{x_0}^x F(-t) dt\big/A^2(t)$  and  $a(-x) = \int_{x_0}^x [1-F(t)] dt\big/ A^{2}(t) +o(a(x))$ as $x\to\infty$ (\cite{Upot}).  This leads to  
   $$\frac1{A(x)}  = -(p-q) \int_{x_0}^x \frac{H(t)}{A^2(t)}dt  + \frac1{A(x_0)} \sim (1-q^{-1}p) a(x)$$  
(with usual interpretation if $p=1/2$),  so that $U_\as(x)V_\ds(x)/x \sim (1-q^{-1}p) a(x)$.
 \end{rem}
 
\begin{lemma}\label{lem4.5} Suppose that  (\ref{stable}) holds with $\al=1$.
 Let  $EX=0$ and  $p\leq 1/2$. Then    
$$\int_1^x\frac{F(-t)dt}{[\ell^*(t)\ell_\sharp(t)]^2} \sim g_\vom(x,x). $$
\end{lemma}
 
\pf   We suppose  $p>0$ for simplicity.  [If $EX=0$ and $p<1/2$, the result follows from (\ref{id_A2}) because of Theorem 1 of  \cite{Uexit2}  that is shown independently of the present issue.]   Recall $[\ell^*(t)]' = P[Z>t] \sim [1-F(t)]/\ell_\sharp(t)$  
 and  then  observe that
\beq
g_\vom(x,x) 
&=& \sum_{k=0}^x v_\ds(k)u_\as   (k) =  \sum_{k=0}^x \frac{v_\ds(k)}{\ell^*(k)}\{1+o_k(1)\}\\
&=&  \frac{V_\ds(x)}{\ell^*(x+1)} \{1+o(1)\}  + \sum_{k=0}^x V_\ds(k)\bigg(\frac1{\ell^*(k)} -\frac1{\ell^*(k+1)}\bigg)\{1+o_k(1)\}\\
&=& \frac1{\tilde A(x)} \{1+o(1)\} + \sum_{k=1}^x \frac{1-F(k)}{[\tilde A(k)]^2}\{1+o_k(1)\},
\eeq
where $\tilde A(t) =\ell^*(t)\ell_\sharp(t)$ and   $o_t(1)\to 0$ as $t\to\infty$. By $P[Z>t]\ell_\sharp(t) \sim 1-F(t)$ 
(see Remark \ref{rem3.1}) and by (\ref{el_el}) it follows that 
\beqn\label{g0}
\frac{1}{\tilde A(x)}= \int_1^x
\frac{F(-t)-1+F(t)}{[\tilde A(t)]^2} dt+ \int_1^x \frac{[1-F(t)]\times o_t(1)}{[\tilde A(t)]^2}dt+ \frac1{\tilde A(1)}
 \eeqn
 and hence, on absorbing $1/\tilde A(1)$ into the second  integral on the RHS, 
\beqn\label{g}
g_{\vom}(x,x) = \int_1^x \frac{F(-t) -1+ F(t)}{[\tilde A(t)]^2}dt \{1+o(1)\} + \int_1^x \frac{H_+(t)}{[\tilde A(t)]^2} \{1+o_t(1)\}dt,
\eeqn
which leads to the equivalences of the lemma. \qed
\v2
For $F$ recurrent,  under (\ref{KMC}) we shall see that $g_\vom(x,x) \sim a(x)$ ($x\to\infty$) in \cite{Uexit2}, so that the equivalence of Lemma \ref{lem4.5}, which we do not use Lemma \ref{lem4.5} in this paper,
gives 
$$\int_1^x\frac{F(-t)dt}{[\ell^*(t)\ell_\sharp(t)]^2}\sim a(x)\sim  \int_1^x\frac{F(-t)dt}{A^2(t)}$$
that  is of interest in the critical case  $\lim a(-x)/a(x) =1$ when we have  little  information about  the behaviour  of $A(x)/ \ell^*(x)\ell_\sharp(x)$. 
 
 \begin{lemma} \label{lem3.7}   If either (C3) or  (C4)  holds, then for any $\de<1$, as  $R\to\infty$
$$1-P_{R-x}(\La_R) =o\bigg(\frac{U_\as(x)}{U_\as(R)}\bigg)  \quad\; \mbox{uniformly for}  \; 0\leq x <\de R.$$
\end{lemma}
\pf
We prove the dual assertion.  Suppose $U_\as$ and $v_\ds$ are s.v.   We show 
\beqn\label{eL3.7}
P_x(\La_R)V_\ds(R)/V_\ds(x)\to 0.
\eeqn  
Since for $\frac13 R \leq x<\de R$, $V_\ds(R)/V_\ds(x)$ is bounded  and $P_x(\La_R)\to 0$ as $R\to\infty$,  we have only to consider  the case  $x<\frac13 R$. 
Putting $R'=\lfloor R/3\rfloor$ and $R''= 2R'$, we have
$$P_x(\La_R) = J_1(R)+J_2(R),$$
where
$$J_1(R)= P_x\big[S_{\sigma(R',\infty)} > R'',\,  \La_{R}\big]\; \;\;  \mbox{and} \quad J_2(R) = P_x\big[S_{\sigma(R',\infty)}\leq R'',\, \La_{R}\big].$$
Obviously 
$$J_1(R) \leq \sum_{z=0}^{R'} g_{\mathbb{Z}\setminus [0,R']}(x,z)P[X > R''-z] 
\leq \big[1-F(R')\big]\sum_{z=0}^{R'} g_{\vom}(x,z).$$
 The last sum is at most    $V_\ds(x)U_\as   (R)$ by Lemma \ref{lem1}.  Using Lemma \ref{lem3.3}(i)
  we accordingly infer that
\[
J_1(R)\,\frac{V_\ds(R)}{V_\ds(x)} \leq V_\ds(R)U_\as   (R) \big[1-F(R') \big] \,\longrightarrow\, 0.
\]
 As for  $J_2$, we decompose $J_2(R) = P_x(\La_{R'}) \sum_{z=R'}^{R''} P_x\big[S_{[R',\infty)}=z\,\big|\, \La_{R'}\big] P_z(\La_R)$. 
 On applying the upper bound (\ref{upp_bd}) to $P_x(\La_{R'})$  it then follows that as $R\to\infty$ 
 \beqn\label{II}
J_2(R)\, \frac{V_\ds(R)}{V_\ds(x)} \leq  \frac{V_\ds(R)}{V_\ds(R')}\sum_{z=R'}^{R''} P_x\big[S_{[R',\infty)}=z\,\big|\, \La_{R'}\big] P_z(\La_R)  \,\longrightarrow\,  0,
\eeqn
for $P_z(\La_R) \to 0$ uniformly for $R'\leq z\leq R''$.  Thus (\ref{eL3.7}) is verified.

 In case (C4),  $V_\ds(x)\sim x^\al/\hat\ell(x)$ and $U_\as   (x)\sim 1/\hat \ell_\sharp(x)$ (see (\ref{5.25}) and (\ref{5.30})) and we can proceed as above.
 \qed
\v2

The next lemma presents some of the results obtained above in a neat form under condition (AS).
\begin{lemma}\label{lem3.8}
 Suppose  (AS) to hold with $\al=1$.
 
   {\rm (i)} \, The following are equivalent
$$\mbox{ (a)\, $P[-\hat Z\geq x]$ is s.v.\;\; \;  (b)\;   $Z$ \;is r.s.  \;\;  (c) \; $\lim P[S_n>0] = 1$,}$$
and each of (a) to (c) above  implies $P[-\hat Z\geq x] \sim \ell_\sharp(x)$ and $V_\ds(x)\sim 1/\ell_\sharp(x)$;   
in particular  these two asymptotic equivalences hold   if  either $EX=0, p<1/2$ or $E|X|=\infty$, $p>1/2$.

 {\rm (ii)} \, If either  of (a) to (c) of {\rm (i)} holds, then  
 \beqn\label{2.10} 
\begin{array}{ll}
\;  \sup_{x>(1+\ep)R}P_x[\sigma_{\{R\}}<T] \to 0\qquad \quad\mbox{for each $\ep>0$; and}\\[2mm]
\left\{\begin{array}{ll}
  P[Z>x] \sim V_\ds(x) [1-F(x)] \quad\quad &\mbox{if} \;\;\; p >0,\\[1mm]
 P[Z>x] = o\big(V_\ds(x)F(-x)\big) \quad  &\mbox{if \, $p=0$}.
\end{array}\right. 
\end{array}
\eeqn
\end{lemma}
\pf  
The equivalence of (a), (b) and  (c) is observed in \cite{Upot}. [The same result can also deduced from what is mentioned in Remark \ref{rem2} combined with our Theorem \ref{thm1}, where one uses the latter to ensure that the slow variation  of $V_\ds$-equivalent to (a)---implies (c).] The other assertion of (i) is contained in Lemma \ref{lem3.1}.

From the slow variation of $V_\ds$ it follows that $\lim_{y\to\infty} P_y[S_{\sigma(-\infty,0]} < -My] =1$ for each  $M>1$ \cite[Theorem 3]{R}, which entails the first relation of (ii).  
The second one  follows immediately from  Lemma \ref{lem3.2}(i) if $p>0$, since one can take $\de=0$ in (\ref{3.1}) under (AS). In case $p=0$, examine the proof of (\ref{3.2}) by noting that  $\sum_yv_\ds(y)[1-F(y)] <\!< \sum_y v_\ds(y)F(-y)$. \qed

\section{Proof of Proposition \ref{prop1.1}}  
\v2
Throughout this section we suppose  (AS) to hold.
Let $p, q$ and $L$ be as in (\ref{stable}). Recall that  $\al\rho$ ranges exactly over the interval $[(\al-1)\vee0,\al\wedge1]$.
 
It is known that  $P[Z>x] $ varies  regularly at infinity with index $- \alpha\rho$ if   $\alpha \rho<1$ (which if $0<\al<2$ is equivalent to $p>0$) and $Z$ is r.s. if  $\alpha \rho=1$  (cf. \cite[Theorem 9]{R} in case $\rho>0$; also Lemma \ref{lem3.1} for the first statement in case $\rho=0$); and in either case
\beqn\label{5.1}
\frac{U_\as   (x)\int_0^xP[ Z>t]dt}{x} \,\longrightarrow \, \frac{1}{\Ga(1+\al\rho)\Ga(2-\al\rho)}
\eeqn
(the dual of (\ref{U/Z})).  We  choose a s.v. function $\ell(x)$  so that
\beqn\label{5.20}
\begin{array}{ll}
 \ell(x) = \ell^*(x) \quad &\mbox{if}\;\; \al\rho=1\\
{\displaystyle P[ Z>x] \sim  \frac{\sin \pi\alpha \rho}{\pi\al\rho } x^{-\alpha \rho}\ell(x)} \quad &\mbox{if} \;\;\al\rho<1.
\end{array}
\eeqn
Here $t^{-1}\sin t$ is understood to equal unity for $t=0$. [Note that in case $\al=1$ with $\rho<1$, $\ell $ is not asymptotically equivalent to $\ell^*$; the constant factor in the second formula above is chosen so as to have the choice  conform to that in (\ref{-1})---see (\ref{5.30}) below.]
As the dual relation we have a s.v.\,function $\hat \ell$ such that
$$ \hat\ell = \hat \ell^*(x) \quad (\al\hat \rho=1); \quad P[-\hat Z>x]/v^0  \sim  \frac{\sin \pi\alpha \hat \rho}{\pi\al\hat\rho} x^{-\alpha \hat\rho} \hat\ell(x) \quad (\al \hat\rho<1).
$$
 Note that  according to  Lemma \ref{lem3.1},
   $$\hat \ell(x) \sim \ell_\sharp(x) \quad \mbox{if} \;\;\rho =1;$$
  similarly  $\ell(x)\sim \hat\ell_\sharp(x)$ if  $\rho=0$ (necessarily $\alpha\leq 1$), where 
  \beqn\label{5.25}
  \hat \ell_\sharp(x) =\al \int_x^\infty \frac{1-F(t)}{t^{1-\al} \hat\ell(t)}dt \qquad (x>0).
  \eeqn
  It then follows that for all $0\leq \rho\leq 1$,
\beqn\label{5.30}
U_\as   (x) \sim  x^{\al\rho}/\ell(x) \quad and \quad V_\ds(x)\sim  x^{\al\hat\rho}/\hat\ell(x).
\eeqn

The two functions $\ell$ and $\hat \ell$ are linked as we are going to see in Lemma \ref{lem6.1} below. 
We bring in the constant
$$\kappa = \frac{\Ga(\al)\pi^{-1}\sin \pi\alpha \rho}{p\Ga(\alpha\rho+1)\Ga(\alpha\hat \rho+1)} =\frac{\Ga(\al)\pi^{-1}\sin \pi\alpha\hat \rho}{q\Ga(\alpha\rho+1)\Ga(\alpha\hat \rho+1)} = \frac{\Ga(\al)[\sin \pi\al\rho  +\sin \pi\alpha\hat \rho]}{\pi\Ga(\alpha\rho+1)\Ga(\alpha\hat \rho+1)},
 $$
where 
only the case $p>0$ or $q>0$ of the first  two expressions above is adopted if $pq=0$;  if $pq\neq 0$ the  two   coincide, namely  
${p}/{q} ={\sin \pi\alpha \rho}/{ \sin \pi\alpha \hat\rho}$ ($pq\neq0$)
as is implicit in the proof of the next result or   directly derived (cf. Appendix (A) of \cite{Upot}).  Note that   $\kappa=0$ if and only if either  $\al=1$ and  $\rho\hat\rho=0$ or $\al=2$. Recall $H(t)$ stands for $P[|X|>t]$.

\begin{lemma} \label{lem6.1}\,  {\rm (i)}  For all $0<\al\leq 2$, \,  $[x^\al H(x)]\big/\ell(x)\hat \ell(x)\to \kappa$, or equivalently
$$U_\as   (x)V_\ds(x)  = \frac{\kappa+o(1)}{H(x)}.$$

{\rm (ii)} If $\al=2$, then  $\ell(x)\hat\ell(x)\sim \int_0^x tH(t)dt$, or what amounts to the same, 
$$U_\as   (x)V_\ds(x)  \sim  \frac{x^2}{\int_0^x tH(t)dt}.$$

{\rm (iii)} If $\al =1$ and $p\neq 1/2$ (entailing $\rho \in \{0,1\}$ so that $\kappa=0$), then
$$U_\as   (x)V_\ds(x)  \sim  \frac{x}{|A(x)|} \sim 
\left\{\begin{array}{ll}
{\displaystyle \frac{x}{|p-q|\int_x^\infty H(t)dt}} \quad\;\; (EX=0),\\[5mm]
{\displaystyle \frac{x}{|p-q|\int_0^x H(t)dt} }\quad\;\; (E|X|=\infty).
\end{array} \right.
$$
 \end{lemma}
 
For the identification of the asymptotic form of  $U_\as   (x)V_\ds(x)$ explicit  in terms of $F$, Lemma \ref{lem6.1} covers all the cases of (AS) other than the case $\al=2p=\rho\vee\hat\rho =1$.

\v2\nin
\pf (iii) follows immediately from Lemmas \ref{lem3.1} and \ref{lem4.4} (see also Remark \ref{rem4.2}). We show (ii) first,  it being used for the proof of (i). 
 
{\sc Proof of} (ii).  Let $\al =2$.  Then $\ell=\ell^*$, $\hat \ell =\hat\ell^*$. Suppose that $EZ=E|\hat Z|= \infty$ so that $\ell^*(x)\wedge \hat\ell^*(x)\to\infty$, otherwise the  proof below merely  being simplified. By (\ref{5.40}) and its dual
\begin{eqnarray}\label{5.00}
(\ell\hat\ell)'(x) &=& \ell^*(x)P[-\hat Z>x]/v^0+ \hat\ell^*(x)P[Z>x]\\
&=& \ell^*(x)\int_0^\infty \frac{F(-x-y)}{\ell^*(y)\vee1}dy\{1+o(1)\} 
+\hat\ell^*(x)\int_1^\infty \frac{1-F(x+y)}{\hat\ell^*(y)}dy\{1+o(1)\}. \nonumber
\end{eqnarray}
We claim that for any constant $\ep \in (0,1)$,
\beqn\label{5.000}
\int_0^x \ell^*(t)dt\int_0^\infty \frac{F(-t-y)}{\ell^*(y)\vee1}dy 
 \left\{
 \begin{array}{ll}
\leq \big[\ep\int_0^x tF(-t)dt + \int_0^x\eta_-(t)dt\{1+o(1)\},\\[2mm]
\geq \int_0^x tF(-t)dt\{1-\ep +o(1)\}.
\end{array}
\right.
\eeqn
By splitting the inner integral of the LHS at $\ep t$ and using the monotonicity of $\ell^*$ as well as of $F(-t)$ it follows that
$$\ell^*(t)\int_0^\infty \frac{F(-t-y)}{\ell^*(y)\vee1}dy \leq \big[\ep tF(-t) +\eta_-(t)\big]\{1+o(1)\},
$$
showing the first inequality of (\ref{5.000}). On the other hand the LHS of (\ref{5.000}) is greater than
\beq
 \int_0^x \ell^*(t)dt\int_0^{x-t} \frac{F(-t-y)}{\ell^*(y)\vee1}dy&=& \int_0^x F(-w)dw\int_0^w\frac{\ell^*(w-y)}{\ell^*(y)\vee1}dy\\
 &\geq& \int_0^x F(-w)dw\int_0^{(1-\ep)w}\{1+o(1)\}dy,
 \eeq
 which yields the second inequality of  (\ref{5.000}). 
 
 Put $L^*(x)= \int_0^x tH(t)dt$. Then (AS) (with $\al=2$) is equivalent to the slow variation of $L^*$ and in this case $m(x)=\int_0^x \eta(t)dt \sim  L^*(x)$ (cf. \cite[Appendix (A)]
{Upot}), which entails $x\eta(x) =o(L^*(x))$. By (\ref{5.000}) and its dual that follows because of symmetric roles of two tails of $F$, integrating the second expression of $(\ell\hat\ell)'(x)$ in (\ref{5.00}) leads to
$$ \ell(x)\hat\ell(x) \bigg\{ \begin{array}{ll} \leq \{1+\ep+ o(1)\} L^*(x),\\[1mm]
 \geq \{1-\ep+ o(1)\} L^*(x), \end{array}$$
which shows the formula of (ii), $\ep$ being arbitrary.
\v2

{\sc Proof of} (i). We have only to consider the case  $0<\rho<1$ (of $0<\al <2$), since for $\rho=1$ the result follows from  Lemma \ref{lem3.3}  and for $\rho=0$ by duality, and  since for $\al =2$  $x^2H(x) = o(\int_0^x tH(t)dt)$ so that the result follows from (ii).
By symmetry we may and do suppose $q>0$ (entailing  $\al\hat\rho<1$).       
   Recalling (\ref{hat-Z}) we have
\beqn\label{5.40}
\frac{P[- \hat Z \geq   x]}{v^0} = \sum_{z=0}^\infty  u_\as   (z) F(-x-z).
\eeqn
We may suppose that $\ell$ and $L_-:= qL$ are differentiable and satisfy $\ell'(t) =o(\ell(t)/t)$ and $L'_-(t) =o(L_-(t)/t)$. 
Summing by parts  deduces that    the sum on the RHS is equal to
\beqn\label{5.50}
\sum_{z=0}^\infty  U_\as   (z) \big[F(-x-z)- F(-x-z-1) \big]. 
\eeqn
On substituting from (\ref{5.30}) and then summing by parts back,  the above sum is  equivalent to 
$$\al\rho\sum_{z=1}^\infty \frac{z^{\alpha \rho-1} L_-(x+z)}{\ell(z)(x+z)^{\alpha}} \sim  \frac{qL(x)}{\ell(x)}\int_0^\infty\frac{\al\rho z^{\alpha\rho-1}}{(x+z)^\alpha}dz.
$$
The  integral  on the RHS equals   $x^{-\al \hat\rho}\alpha\rho B(\al\rho,\alpha\hat \rho)$, where $B(t,s)=\Ga(t)\Ga(s)/\Ga(t+s)$ (the Bessel function). We accordingly  obtain
\beqn\label{hat_Z2}
P[-\hat Z\geq x]/v^0 \sim q\al\rho B(\alpha\rho,\alpha\hat \rho) x^{-\al\hat\rho\,}L(x)/\ell(x).
\eeqn
Thus by (\ref{U/Z}) (the dual  of  (\ref{5.1})) and the identity $\Ga(1-t)\Ga(1+t) =\pi t/\sin t$ ($|t|<1$)
\beqn\label{5.70}
V_\ds(x) \sim \frac{\sin \pi\alpha \hat\rho}{(\pi \alpha^2\rho \hat\rho )B(\alpha\rho,\alpha\hat \rho)}\cdot\frac{x^{\alpha\hat \rho}\ell(x)}{qL(x)} 
=\frac{\Ga(\al)\pi^{-1}\sin \pi\alpha \hat\rho}{\Ga(\alpha\rho+1)\Ga(\alpha\hat \rho+1)}\cdot\frac{x^{\alpha\hat \rho}\ell(x)}{qL(x)},
\eeqn
which combined with $U_\as   (x) \sim x^{\al \rho}/ \ell(x)$ shows the asserted convergence. 

Proof of Lemma \ref{lem6.1} is complete. \qed
\v2
\begin{rem}\label{rem5.1} Under  (AS) the equivalence (\ref{hat_Z2}) holds  unless  $q=\hat\rho=0$ with the understanding that $sB(s,t) =1$ for $s=0, t>0$ and $B(s,0)=\infty$ for $s>0$,  and $\sim$ is interpreted in the usual way if the constant factor of its right member equals zero or infinity. The proof is carried out in the 
following three cases
\vskip1mm
(a) either $\alpha=\rho =1>p$ or $\alpha=\hat\rho=1$;\quad  (b) $\alpha<1=q$;\quad (c) $\alpha>1= p$.
\vskip1mm\noindent
Suppose  (a) holds with $\rho=1$. Then $\ell =\ell^*$ and $P[-\hat Z>x]/v^0\sim \ell_\sharp(x) >\!> x\mu_-(x)/\ell^*(x)$. Since $q>0=\hat\rho$, this shows (\ref{hat_Z2}).  If  (a) holds with $\hat\rho=1$, then $\ell=\hat\ell_\sharp$ and $U_\as(x)\sim 1/\hat\ell_\sharp(x)$, and from the representation
$
P[-\hat Z>x]/v^0 = \sum_{y=0}^\infty u_\as(y)\mu_-(x+y)
$
we deduce that
$$P[-\hat Z>x]/v^0 =\bigg\{
\begin{array}{ll} U_\as(x)\mu_-(x)\{1+o(1)\}  \sim qx^{-1}L(x)/\ell(x) \;\; &\mbox{if}\; q>0,\\
o(\mu(x)/\ell(x)) &\mbox{if}\; q=0,
\end{array}
$$
showing (\ref{hat_Z2}).
In case (b) we also have $\ell=\hat\ell_\sharp$ and $U_\as(x)\sim 1/\hat\ell_\sharp(x)$ and obtain  (\ref{hat_Z2}) in the same way as above. As for  the case (c) we  have  $\alpha\rho=1$, $\ell=\ell^*$, and $u_\as(x)\sim 1/\ell^*(x)$, and see  that  $\sum_{y=0}^\infty u_\as(y)\mu_-(x+y) <\!< U_\as(x)\mu(x) \sim x^{-\alpha\hat\rho}L(x)_\ell(x)$ since
 $\mu_-/\mu \to0$.
\end{rem}

\v2
The asymptotic form of $V_\ds U_\as$  given in Lemma \ref{lem6.1} may also be expressed as 
\beqn\label{nUV}
\lim U_\as   (c_n)V_\ds(c_n)/n \,  \longrightarrow\,   \mbox{$\kappa$\, or\,  $1/2$ \quad according as\, $\al <2$\, or\, $\al =2$}
\eeqn
(if $\al =2$ the sequence ($c_n$) is determined by $c^2_n/L(c_n) \sim  n$, where $L(x)= \int_{-x}^x t^2dF(t)$) and known if $0<\rho <1$ \cite[Eq(15)]{D_lb}, \cite[Eq(15, 31)]{VW} apart from the explicit expression of the limit value.

\begin{lemma}\label{lem6.2}  Suppose $\rho>0$. Then for any $\ep>0$,  there exists a constant  $M>0$ 
such that  for all sufficiently large $R$,
 $$E_x\big[V_\ds(S_{\sigma(R,\infty)}); Z(R) \geq  MR,\, \La_R\big] \leq \ep V_\ds(x)
\quad \; (0\leq x<R),$$
and,  if $\al\rho=1$, 
$$\lim_{R\to\infty} \sup_{0\leq x<R} \,  \frac1{V_\ds(x)}E_x\big[V_\ds(S_{\sigma(R,\infty)}); Z(R)\geq \ep R,\La_R\big] =0.$$
 \end{lemma}
 \pf \, As in the proof of Lemma \ref{lem3.4} but with  $V_\ds(x) \sim x^{\al(1-\rho)}/\hat\ell(x)$ in place of $V_\ds(x)\sim1/\ell_\sharp(x)$,
 one deduces that if $\rho>0$, for all sufficiently large $R$,
 \beqn\label{5.11}
\sum_{w\geq R+MR} p(w-z) V_\ds(w) \leq C_*\frac{V_\ds(R)H(R)}{M^{\al\rho}}\quad\; (0\leq x<R)
 \eeqn
 with some constant  $C_*$ for which one can take $2/\alpha\rho$. Using  Lemmas \ref{lem1} and \ref{lem6.1}(i) in turn one then obtains that for any $M>0$
\begin{eqnarray}\label{5.12}
&&\frac1{V_\ds(x)}E_x\big[V_\ds(S_{\sigma(R,\infty)}); Z(R)  \geq M R,\La_R\big]   \\
&&\qquad= \frac1{V_\ds(x)}\sum_{z=0}^{R} \sum_{w\geq (M+1)R} g_{B(R)}(x,z)V_\ds(w) p(w-z)\nonumber\\
&&\qquad\leq
C_*\frac{U_\as   (R)V_\ds(R)\mu(R)}{M^{\al\rho}} \nonumber\\
&&\qquad = C_*\frac{\kappa+o(1)}{M^{\al\rho}}, \nonumber
\end{eqnarray}
 so that the first half of the lemma is ensured by taking
 \beqn\label{5.13}
 M=[2C_*\kappa/\ep]^{1/(\al\rho)}.
 \eeqn
 If  $\al\rho=1$, we have either $\kappa=0$ or $1<\al<2$  with $p=0$ and   the second half follows from (\ref{5.12}),  for if $p=0$,  $H(x)$ can be replaced by $o(H(x))$.
 \qed
 
 \v2
\begin{lemma}\label{lem6.3}\,  Suppose $ \rho >0$. Then for some constant $\theta>0$,  $P_x(\La_R) \geq \theta V_\ds(x)/V_\ds(R)$ for $0\leq x< R$ and for $R$ large enough.
If $\al \leq 1$, then  $\theta$ can be chosen so as to approach unity as $\rho\to 1$, and for any $\de >0$, for $x\geq \de R$,
$$P_x(\La_R) \to 1\quad \mbox{as\; $R\to\infty$\; and\; $\rho \to 1$\; in this order}.
$$
 \end{lemma}
 \pf 
  Because of the identity (\ref{mart}) (due to the optional stopping theorem) Lemma \ref{lem6.2} shows that for any $\ep>0$ we can choose a constant $M>0$ so that $$
(1-\ep) V_\ds(x) \leq E_x\big[V_\ds(S_{\sigma(R,\infty)}); Z(R) < MR\,\big|\, \La_R\big] P_x(\La_R).$$ 
Since the conditional expectation 
 on the RHS is less than  $V_\ds(R+MR) \sim (1+M)^{\al \hat\rho}V_\ds(R)$
 we have the lower bound of $P_x(\La_R)$ with $\theta = (1-\ep)(1+M)^{-\al\hat\rho}$. Noting that $M$ may be determined by (\ref{5.13}) and that as  $\rho \to 1$ (possible only when  $\al\leq 1$), 
 $\theta\to1-\ep$, hence $\theta \to 1$ for $\ep$ may be arbitrarily small.  Since one also has  
 $V_\ds(x)/V_\ds(R) \sim (x/R)^{\al\hat\rho} \to 1$ (uniformly for $\de R<x<R$) one can concludes the second half of the lemma.  \qed

\begin{lemma} \label{lem6.4}\, Suppose $0<(\al\vee1)\rho<1$. Then  for any $\de<1$, uniformly for $0\leq x < \de R$,  as $R\to\infty$ and $\ep \downarrow 0$ (interchangeably) 
$$P_x\big[Z(R)  \leq \ep R\,\big|\, \La_R\big]\to 0.$$
\end{lemma}
\pf  By the restriction as to $(\al\vee1) \rho$ it follows that $\al<2$ and
\beqn\label{5.14}
\begin{array}{ll}
{\rm (a)}\;\;  p>0\; \mbox{so that $1-F$ is regularly varying;}\\
{\rm (b)}\;\; P[Z>x] \; \mbox{is regularly varying with index $\in (0,1)$;}\\
{\rm (c)}\;\; V_\ds(x)/V_\ds(R)\leq CP_x(\La_R)\; \mbox{and for each $\ep>0$, 
$\inf_{R\geq1}\inf_{x\geq \ep R} P_x(\La_R)>0$.} \qquad
\end{array}
\eeqn
The reasoning for (a) is already mentioned, while (b) and (c) are due to (\ref{5.20}) and Lemma \ref{lem6.3}, respectively.  According to \cite[Section XIV.3]{F} (b) implies that $P[Z(r)<\ep r]\to0$ as $\ep\to 0$ uniformly for $r\geq 1$. This in particular shows the asserted convergence of the lemma to be valid for $x\geq \frac14 R$, therein  $P_x(\La_R)$ being bounded away from zero by (c).

For $x <\frac14 R$, take a constant $0<h <1/4$ arbitrarily and let ${\cal E}_{\ep, R}$ stand for the event  $\{(1-h)R\leq S_{\sigma[R/4,\infty)} <(1+\ep)R\}$. It holds that
\beqn\label{5.15}
\La_R\cap \{Z(R)<\ep R\} \subset \La_{R/4}\cap \big[ \{S_{\sigma[R/4,\infty)} <(1-h)R, Z(R)<\ep R\}\cup {\cal E}_{\ep, R}\big],
\eeqn
and that
\beqn\label{5.16}
P_x(\La_{R/4} \cap  {\cal E}_{\ep, R}) \leq 
 \sum_{0\leq y <R/4} g_{\vom}(x,y)P[(1-h)R-y\leq X< (1+\ep) R-y].
\eeqn
Put
$$\om_{R}(h,\ep) = \sup_{y<R/4} \frac{P[(1-h)R-y\leq X< (1+\ep) R-y]}{P[X\geq (1-h)R-y]}.
$$
One sees that $\om_{R}(h,\ep) \leq \sup_{z<R/2} P[z\leq X < z+ (h +\ep)z]/P[X\geq z]$, hence the regular variation of  $1-F$ entails that $\om_{R}(h,\ep) \to 0$ as $h+\ep\to0$ uniformly for $R\geq 1$. Since in 
(\ref{5.16}) the probability under the summation sign is less than $\om_{R}(h,\ep)[1-F(R/2)]$, by
Lemmas \ref{lem1} and \ref{5.1} the sum is evaluated to be at most a constant multiple of
$$
\om_{R}(h,\ep) \big[1-F({\textstyle \frac12}R)\big] V_\ds(x)U_\as   ({\textstyle \frac12}R)\leq \om_{R}(h,\ep) \frac{V_\ds(x)}{ V_\ds(\frac12 R)} \{\kappa+o(1)\}.
$$
By (\ref{5.14}c) and regular variation of $V_\ds$,  $V_\ds(x)/V_\ds(\frac12 R)\leq CP_x(\La_R)$. Hence 
\beqn\label{5.17}
P_x[\La_{R/4}\cap {\cal E}_{\ep, R}\,\big|\, \La_R\big] \leq C_2\, \om_{R}(h,\ep)    \quad\mbox{for}\quad  0\leq x< R/4.
\eeqn
Now by (\ref{5.15}),   $P_x[Z(R) \leq \ep R\,|\, \La_R] $ is less than 
\[
P_x[\La_{R/4}\cap {\cal E}_{\ep, R}\,\big|\, \La_{R}\big]
+ \sup_{R/4\leq y< (1-h)R}P_y[Z(R)<\ep R]\frac{P_x(\La_{R/4})}{P_x(\La_R)}.
\]
By (\ref{5.14}c) again, $P_x(\La_{R/4})/P_x(\La_R) =1/P_x[\La_R\,|\,\La_{R/4}] <C$. By (\ref{5.17}),   for any $\ep_0>0$,  one  can choose $h>0$ so that the first term is less than $\ep_0$ for all $0<\ep<h$ and all sufficiently large $R$, while for $h$ thus chosen the second term is less than $\ep_0$ for $\ep$ small enough and all $R$ large enough.
Thus the convergence asserted in the lemma follows. \qed

\begin{lemma} \label{lem6.5}  Let $\al \leq 1$. For any $\de<1$,
$$\sup_{0\leq x <\de R} \frac{P_x(\La_R)V_\ds(R)}{V_\ds(x)} \to 0 \quad\; \mbox{as\; $R\to\infty$\; and\; $\rho\to0$\; in this order;}$$
if $\rho=0$, the above supremum approaches zero as $R\to\infty$.
\end{lemma}
\pf  The case  $\rho=0$ is a special case  of Lemma \ref{lem3.7}. 
 The proof of the rest  is made by examining that of Lemma \ref{lem3.7}.  Let $\al\leq 1$. Suppose  $0<\rho<1$. 
Since for $\frac13 R \leq x<\de R$, $V_\ds(x)/V_\ds(R)$ is bounded away from zero and $P_x(\La_R)\to 0$ as $R\to\infty$, $\rho\to 0$ in this order because of the dual of the last assertion of Lemma \ref{lem6.3} we have only to consider  the case  $x<\frac13 R$. Putting  $J_1(R)$ and $J_2(R)$ as in the proof of Lemma \ref{lem3.7} so that $P_x(\La_R) = J_1+J_2$.
 Using Lemmas \ref{lem1} and  \ref{lem6.1}(i), 
  we have 
\beqn\label{5.191}
J_1(R)\,\frac{V_\ds(R)}{V_\ds(x)} \leq V_\ds(R)U_\as   (R) \big[1-F(R') \big] \leq 4p\kappa
\eeqn
for $R$ large enough. As before we  see  that as $R\to\infty$ and $\rho\to 0$ in turn,
\[
J_2(R), \frac{V_\ds(R)}{V_\ds(x)} \leq  \frac{V_\ds(R)}{V_\ds(R')}\sum_{z=R'}^{R''} P_x\big[S_{[R',\infty)}=z\,\big|\, \La_{R'}\big] P_z(\La_R)\to 0,
\]
for, in this limit, $P_z(\La_R) \to 0$ uniformly for $R'\leq z\leq R''$. Since $p\kappa \to 0$ as $\rho\to 0$, this together with (\ref{5.191}) shows the first half of the lemma. 
 \qed

\v2\v2
{\bf  Proof of   Proposition \ref{prop1.1}.} Let $0< (\al\vee1) \rho<1$.  The  lower bound of $P_x(\La_R)$ asserted in (i) follows from Lemma \ref{lem6.3}. For proof of the upper bound let $1/2 <\de <1$.  Lemma \ref{lem6.4} entails 
\beqn\label{pr_(i)}
P_x\big[Z(R) > \ep R\,\big|\, \La_R\big] \geq 1/2\qquad  (x<\de R)
\eeqn  
for some $\ep>0$. 
Writing the equality (\ref{mart}) as
$$\frac{V_\ds(x)}{P_x(\La_R)} =\sum_{y=R+1}^\infty  P_x[S_{N(R)}=y\,|\,
 \La_R] V_\ds(y),$$  
one deduces from   (\ref{pr_(i)}) that  the sum above is larger than $[V_\ds(R) + V_\ds((1+\ep)R)]/2$. By (\ref{5.30})  
$\limsup V_\ds((1+\ep) R)/V_\ds(R) =(1+ \ep)^{\al\hat\rho}>1$, for  $(\al\vee1) \rho<1$ entails $\hat \rho >0$. Hence $P_x(\La_R) \leq \big(2/[1+(1+\ep)^{\al\hat\rho}]\big)V_\ds(x)/V_\ds(R)\{1+o(1)\}$.
 This verifies the upper bound.

(ii) follows from Lemmas \ref{lem6.3} and \ref{lem6.5}.
\qed

\section{Appendix}

\;\; (A)  Here we present three lemmas involving s.v. functions, of which the first, second and third ones are used in the proofs of Lemmas \ref{lem3.1}, \ref{lem3.2} and \ref{lem3.6}, respectively. 
\begin{lemma} \label{lem8.1} Let $\al$ be a positive constant,  $\ell(t), t\geq 0$ a s.v.   function,
  $u(t)$ and $G(t)$  positive measurable functions of  $t\geq 0$.  Suppose  that   $\int_0^x u(t)dt  \sim x^\al/\al \ell(x)$ ($x\to\infty$),  $G(t)$ is  non-increasing  and both $u$ and $1/\ell$ are locally bounded, and  put
$$\tilde u(t) =\frac{t^{\al-1}}{\ell(t)},\quad  h_\sharp(x)=\int_x^\infty \tilde u(t)G(t)dt,$$
$$ h(x)=\int_0^\infty u(t)G(x+t)dt, \quad  \mbox{and}\quad \tilde h(x)=\int_0^\infty \tilde u(t)G(x+t)dt,$$
If one of $h(1), \tilde h(1)$ and $h_\sharp(1)$ is finite, then so are the other two and
\vskip1mm
{\rm (i)} \quad if either $h$ or $\tilde h$ is regularly varying with index $>-1$, then  $h(x)\sim \tilde h(x);$
\vskip1mm
{\rm (ii)}   \quad if either $\tilde h$ or  $h_\sharp$ is s.v., then  $h(x)\sim h_\sharp(x)$.
\end{lemma}
\nin
\pf\, We may suppose $\ell$ is positive and continuous. Integrating by parts verifies that $\int_0^x \tilde u(t)G(t) \asymp \int_0^x u(t)G(t)dt$ so that  $h(1), \tilde h(1)$ and $h_\sharp(1)$ are finite if either of them is finite.

Note that $\tilde h(1)<\infty$ entails $G(t)=o(t^{-\al}\ell(t))$, so that $G(t)\int_0^tu(s)ds \to0$ ($(t\to\infty$). Put $U(x) =\int_0^t u(s)ds$ and $\tilde U(t)= \int_0^t \tilde u(s)ds$. Then, interchanging the order of integration and  integrating by parts in turn one deduces 
$$\int_0^x h(t)dt  =\int_0^\infty u(t)dt \int_t^{x+t} G(s)ds = \int_0^\infty [G(t)-G(x+t)]U(t) dt,$$
and similarly
$$\int_0^x \tilde h(t)dt =  \int_0^\infty [G(t)-G(x+t)] \tilde U(t)dt.$$
Since $G$ is monotone and $U(t)\sim \tilde U(t)$, these identities together show
 $$\int_0^x h(s)ds \sim \int_0^x \tilde h(s)ds,$$
 provided $\int_0^x h(s)ds$ or $\int_0^x \tilde h(s)ds$ diverges to infinity. It therefore follows that if  one of  $h$ or $\tilde h$ is regularly varying with index $>-1$, then so is the other and  $ \tilde  h(x)\sim h(x)$ owing to 
 the monotone density theorem \cite{BGT}. Thus (i) is verified.

 For proof of (ii), suppose either $\tilde h$ or $h_\sharp$ is s.v.   Take  $M>2$ large. Then we see 
$$\tilde h(x) \geq \int_0^{(M-1) x} \tilde u(t)G(x+t)dt = (M-1)^\al\frac{x^\al}{\al \ell(x)} G(Mx)\{1+o(1)\}
$$
and
\beqn\label{uG/Gh}
\tilde h(Mx) = \bigg(\int_0^{x} +\int_x^\infty\bigg)\tilde u(t)G(Mx+t)dt \leq \frac{x^\al}{\al \ell(x)} G(Mx) +
h_\sharp(x).
\eeqn
If $\tilde h$ is s.v. so that  $\tilde h(Mx)\sim \tilde h(x)$, these inequalities lead to  
$h_\sharp(x)\geq \tilde h(x)\big\{1+O(M^{-\al})+o(1)\big\}$.  On the other hand if $h_\sharp$ is s.v., then  
$$h_\sharp(x)\sim h_\sharp(x/2) \geq \int_{x/2}^x \tilde u(t)G(t)dt 
\geq (1-2^{-\alpha}) \frac{G(x)x^{\alpha}}{\alpha \ell(x)},$$
so that by (\ref{uG/Gh}) applied wth  $x/M$ in place of  $x$, one obtains 
$\tilde h(x) \leq \big\{1+  O(M^{-\alpha}) +o(1)\big\}h_\sharp (x)$. Since $M$ can be made arbitrarily large, one concludes $\tilde h(x) \leq h_\sharp(x)\{1+o(1)\}$ in either case.
The reverse inequality is verified as follows. If $\al\geq 1$, then taking $\ep\in (0,1)$ so small that for
$0<s<1$, $(1-\ep  s)^{\al-1}> 1-\al\ep$, one deduces
$$\tilde h(\ep x) \geq \int_x^\infty \frac{(t-\ep x)^{\al-1}}{\ell(t-\ep x)} G(t)dt \geq (1-\al\ep)h_\sharp(x)\{1+o(1)\},
$$
which shows  $\tilde h(x)\geq h_\sharp(x)\{1+o(1)\}$. 
If $\al <1$, one may suppose $\ell$ to be normalised so that for all sufficiently large $x$, $\tilde u(t)\geq \tilde u(x+t)$ ($t\geq0$),  hence $\tilde h(x) \geq h_\sharp(x)$ as well.  
Now we can conclude that  $\tilde h(x)\sim h_\sharp(x)$,  and hence  $h(x)\sim h_\sharp(x)$ by virtue of  (i).
\qed

\begin{lemma}\label{lem8.2} Let $V(x)$ and $\mu(x)$ be positive  functions on $x\geq 0$. If $\mu$ is non-increasing, $\limsup \mu(\lambda x)/\mu(x) <1$ for some  $\lambda >1$ and $V$ is non-decreasing and s.v., then 
$\int_1^\infty \mu(t)V(t)<\infty$  and 
$$\int_x^\infty \mu(t)dV(t) = o\big(V(x)\mu(x)\big)\quad (x\to\infty).$$
\end{lemma}
\pf\, We may suppose that $V$ is continuous and there exists constants $\de<1$ and $x_0$ such that $\mu(\lambda x)<\de \mu(x)$ for $x\geq x_0$.  Let $x>x_0$. It then follows that
\beqn\label{7.1}
\int_x^\infty \mu(t)dV(t) =\int_{n=0}^\infty \int_{\lambda^nx}^{\lambda^{n+1}x} \mu(t)dV(t) \leq \mu(x)\sum_{n=0}^\infty \de^n [V(\lambda^{n+1}x) - V(\lambda^nx)].
\eeqn
Take $\ep >0$ so that $(1+\ep)\de <1$. Then for all sufficiently large $x$, $V(\lambda x)/V(x)\leq 1+\ep$ so that $V(\lambda^n x)/V(x) <(1+\ep)^n$, hence  $\sum_{n=0}^\infty \de^n V(\lambda^n x) \leq CV(x)$, which shows that the last member in (\ref{7.1}) is $o(\mu(x)V(x))$, for $V(\lambda^{n+1} x) - V(\lambda^2 x)= o((V(\lambda^n x))$  uniformly in $n$. \qed

\begin{lemma}\label{lem8.3}
 Let $v(t)$ and $G(t)$ be non-negative measurable functions on $t\geq 0$ such that $V(t):=\int_0^t v(s)ds$ is s.v. at infinity, $G$ is non-increasing and $\int_1^\infty v(t)G(t)dt<\infty$. Then as  $t\to\infty$
 $$\frac{\int_t^\infty v(s)G(s)ds}{V(t)G(t) +\int_t^\infty V(s)G(s)ds/s} \,\longrightarrow\, 0.$$
 \end{lemma}
 \pf
  Let $\tilde V$ be a normalized version of $V$: $\tilde V(t)\sim V(t)$, $ V$ is differentiable and  
  $$\tilde v(t):= \tilde V'(t)= o(V(t)/t).$$ We show that
  \beqn\label{8.2}
  \int_t^\infty v(s)G(s)ds =\int_t^\infty \tilde v(s)G(s) ds\{1+o(1)\} + o\big(V(t)G(t)\big),
  \eeqn
which implies the asserted result. Put 
$D(t) =\int_t^\infty v(s)G(s)ds,\,  E(t)=-\int_t^\infty V(s)G(s)ds.$ 
Then $E(t)= V(t)G(t)+D(t)$.
  Since $E(t)\sim -\int_t^\infty \tilde V(s)dG(s)$, this entails
\beq
D(t) &=& -\tilde V(t)G(t)\{1+o(1)\} -\int_t^\infty \tilde V(s)dG(s) +o(E(t))\\
&=& \int_t^\infty \tilde v(s)G(s)ds + o\big(V(t)G(t)\big) ;+ o(E(t)).
\eeq
Since $o(E(t)) = o(V(t)G(t)) +o(D(t))$, we obtain (\ref{8.2}).
\vskip2mm

(B) \; Let $F$ be transient so that we have the Green kernel 
$G(x):= \sum_{n=0}^\infty P[S_n=x]<\infty$.
For  $y\geq 0, x \in \mathbb{Z}$,
 \beqn\label{8.3}
 G(y-x) -g_{\vom}(x,y) = \sum_{w=1}^\infty P_x[S_T=-w]G(y+w).
 \eeqn
 According to the Feller-Orey renewal theorem \cite[Section XI.9]{F}, $\lim_{|x| \to \infty} G(x)=0$ (under $E|X|=\infty$),  showing that the RHS above tends to zero as $y\to\infty$ (uniformly in $x\in \mathbb{Z}$), in particular $\lim g_{\vom}(x,x) =G(0) =1/P[\sigma_{\{0\}}=\infty]$. It also follows that  $P_x[\sigma_{\{0\}}<\infty]= G(-x)/G(0) \to0$.


\vskip4mm

\begin{thebibliography}{99}

\setlength{\baselineskip}{9pt}


\bibitem{Bt} J. Bertoin, L\' evy Processes, Cambridge Univ. Press, Cambridge  (1996).





\bibitem {BGT} N.H. Bingham, G.M. Goldie and J.L. Teugels, Regular variation,  Cambridge Univ. Press,  Cambridge, 1989.




  






\bibitem {D_lb} R.A. Doney, Local behaviour of first passage probabilities, Probab. Theor. Rel. Fields, {\bf 152}, (2012), 559-588. 


\bibitem {Dsp}  R.A. Doney, Conditional limit theorems for asymptotically stable random walks, Z. Wahrsch.
 Verw. Gebiete {\bf 70} (1985), 351-360. 
















\bibitem {E} K. B. Erickson, The strong law of large numbers when the mean is undefined, TAMS. {\bf 185} (1973), 371-381.

\bibitem{F} W. Feller, An Introduction to Probability Theory and Its Applications, vol. 2, 2nd edn.  John Wiley and  Sons, Inc. NY. (1971)


\bibitem{GM}  P. Griffin and T. McConnell, Gambler's ruin and the first exit position of random walk from large spheres, Ann. Probab. {\bf 22} (1994), 1429-1472.











\bibitem {K} H. Kesten,  Random walks with absorbing barriers and Toeplitz forms, Illinois J. Math. {\bf 5} (1961), 267-290. 



 
 \bibitem {KM/ifp} H. Kesten and R. A. Maller, Infinite limits and infinite limit points of random walks and trimmed sums,  Ann. Probab., \textbf{ 22} (1994), 1473-1513.
 
  \bibitem {KMdivg} H. Kesten and R. A. Maller, Divergence of a random walk through deterministic and random subsequences, J. Theor. Probab., \textbf{10} (1997), 395-427.

\bibitem {KM0} H. Kesten and R. A. Maller, Stability and other limit laws for exit times of random walks from a strip or a half line, Ann. Inst. Henri Poincar\'e, \textbf{ 35} (1999), 685-734.

 
 
 

 \bibitem {P}  S. C. Port, The exit distribution of an interval for completely asymmetric stable processes, 
 Ann. Math. Statist. {\bf 41} (1970), 39-43.

 
 \bibitem {R} B.A. Rogozin,  On the distribution  of the first ladder moment and height and fluctuations of a random walk, Theory Probab. Appl. \textbf{ 16} (1971), 575-595.
 
 \bibitem {R2}  B. A. Rogozin. The distribution of the first hit for stable and asymptotically stable walks on an interval (in Russian). Theor. Probab. Appl. {\bf 17}, 342 – 349 (1972).
 

 





\bibitem {S} F. Spitzer, Principles of Random Walks, Van Nostrand, Princeton, 1964. 




\bibitem{U1dm}  K. Uchiyama,  One dimensional lattice random walks with absorption  at a point / on a half line.  J. Math. Soc. Japan, {\bf 63} (2011), 675-713.




 
      
 
 \bibitem {Uexp}  K. Uchiyama, On the ladder heights of random walks attracted to stable laws of exponent 1, Electron. Commun. Probab. {\bf 23} (2018), no. 23, 1-12.  doi.org/10.1214/18-ECP122
 

\bibitem {Urenwl} K. Uchiyama, A renewal theorem for relatively stable variables.  Bull. L. Math. Society. 
 {\bf 52} (2020), 1174-1190. 

   \bibitem{Upot} K. Uchiyama,  Estimates of Potential functions of 
  random walks on $\mathbb{Z}$ with zero mean and infinite variance and their applications.
  (preprint available at: http://arxiv.org/abs/1802.09832.)

 \bibitem {Uladd}  K. Uchiyama,  The potential function and ladder variables  of a recurrent random walk on $\mathbb{Z}$  with infinite variance. Electron. J. Probab. {\bf 25} (2020)

 \bibitem {Uexit2}  K. Uchiyama,  The two-sided exit problem  for a random walk on  $\mathbb{Z}$
 with infinite variance II. (preprint available at: http://arxiv.org/abs/2102.04102.)

\bibitem{VW}  V. A. Vatutin and V. Wachtel, Local probabilities for random walks conditioned to stay positive, Probab. Theory Rel. Fields, {\bf 143} (2009), 177-217





\end{thebibliography}
\end{document}